\documentclass{article}
\usepackage{amsmath}[1996/11/01]
\usepackage{amssymb,amsthm,amsxtra}
\usepackage{epsfig}

\newtheorem{thm}{Theorem}
\newtheorem{cor}[thm]{Corollary}
\newtheorem{lem}[thm]{Lemma}
\newtheorem{prop}[thm]{Proposition}
\theoremstyle{remark}
\newtheorem{rem}{Remark}

\theoremstyle{definition}

\numberwithin{equation}{section}
\numberwithin{thm}{section}

\DeclareMathOperator{\Exp}{\mathbb E}

\DeclareMathOperator{\Prob}{\mathbb P}

\newcommand{\C}{\mathbb C}
\newcommand{\R}{\mathbb R}
\newcommand{\Z}{\mathbb Z}
\newcommand{\N}{\mathbb N}

\begin{document}

\title{{\bf Optimal tail estimates 
for directed last passage \\ site percolation 
with geometric random variables}}
\author{{\bf Jinho Baik}\footnote{
Deparment of Mathematics,
Princeton University, Princeton, NJ 08544, USA, 
jbaik@math.princeton.edu}
\footnote{Institute for Advanced Study,
Princeton, NJ 08540, USA},
 \ \ 
{\bf Percy Deift}\footnote{Courant Institute of Mathematical Sciences, 
New York University, New York, NY 10012, USA, deift@cims.nyu.edu},
 \ \
{\bf Ken McLaughlin}\footnote{Department of Mathematics, 
University of North Carolina, Chapel Hill, NC 27599, USA, mcl@amath.unc.edu}\\
 \ \
{\bf Peter Miller}\footnote{Department of Mathematics, 
University of Michigan, Ann Arbor, MI 48109, USA,  millerpd@umich.edu},
 \ \
and \ \ {\bf Xin Zhou}\footnote{Department of Mathematics, Duke University, 
Durham, NC 27708, USA, zhou@math.duke.edu}}

\maketitle

\begin{abstract}

In this paper, we obtain optimal uniform lower tail estimates for 
the probability distribution of the properly scaled length of 
the longest up/right path of the last passage site percolation model 
considered by Johansson in \cite{kurtj:shape}.
The estimates are used to prove a lower tail moderate deviation result 
for the model.
The estimates also imply the convergence of moments, 
and also provide a verification of the universal scaling law relating 
the longitudinal and the transversal fluctuations of the model.

\end{abstract}


\section{Introduction}

In \cite{kurtj:shape}, Johansson considered directed 
last passage site percolation on $\Z_+^2=\{ (m,n): m,n \in\N\}$ 
with geometric random variables.
More precisely, for $(i,j)\in \Z_+^2$, let 
$w(i,j)$ be independent, identically distributed geometric random variables
with 
\begin{equation}
  \Prob( w(i,j)=k) = (1-t^2)(t^2)^k, \qquad k=0,1,2,\cdots,
\end{equation}
and $0<t<1$.
An up/right path $\pi$ from $(1,1)$ to $(M,N)$ is, by definition,
a collection of sites $\{(i_k, j_k): k=1,2,\cdots, n\}$, $n:=M+N-1$
such that $(i_1, j_1)= (1,1)$, $(i_n, j_n)=(M,N)$ and
$(i_{k+1}, j_{k+1})-(i_k, j_k)= (1,0)$ or $(0,1)$.
Let $(1,1)\nearrow (M,N)$ be the (finite) set of all such up/right paths
from $(1,1)$ to $(M,N)$.
Now define the maximal `length',
\begin{equation}
  G(M,N) := \max \bigl\{ \sum_{(i,j)\in \pi} w(i,j):
\pi\in (1,1)\nearrow(M,N)\bigr\}.
\end{equation}

Fix $\gamma\ge 1$ and set $M=[\gamma N]$, the integer part of $\gamma N$.
The main result in \cite{kurtj:shape} is the following:
for any $x\in\R$, 
\begin{equation}\label{e-CLT}
  \lim_{N\to\infty} \Prob\biggl(
\frac{ G([\gamma N], N) - \frac1{a_0}N}{b_0N^{1/3}} \le x\biggr)
= F(x),
\end{equation}
where 
\begin{equation}\label{e-a}
  a_0 = \frac{1-t^2}{t((\gamma+1)t+2\sqrt{\gamma})},
\end{equation}
and 
\begin{equation}\label{e-b0}
  b_0 =
\frac{t^{1/3}\gamma^{-1/6}}{1-t^2}
(t+\sqrt\gamma)^{2/3}(1+t\sqrt\gamma)^{2/3},
\end{equation}
and where $F(x)$ is the Tracy-Widom distribution \cite{TW1} for the largest
eigenvalue of a random matrix chosen from the Gaussian unitary ensemble (GUE).
In addition to \eqref{e-CLT}, Johansson also proved 
large deviation results, 
\begin{eqnarray}
\label{e-ld1}
  \lim_{N\to\infty} \frac1{N^2} \log  \Prob\bigl(
G([\gamma N], N) \le N(\frac1{a_0}-y) \bigr) = -\ell(y), \\
\label{e-ld2}
  \lim_{N\to\infty} \frac1{N} \log  \Prob\bigl(
G([\gamma N], N) \ge N(\frac1{a_0}+y) \bigr) = -i(y), 
\end{eqnarray}
for some explicit positive functions $\ell(y)$ and $i(y)$, 
$y >0$.

The result \eqref{e-CLT} parallels an earlier result in \cite{BDJ} 
on the length of the longest increasing subsequence $\ell_N(\sigma)$ 
of a random permutation $\sigma$ of $N$ letters.
The main result in \cite{BDJ} is the following:
for any $x\in\R$, 
\begin{equation}\label{e-BDJ}
  \lim_{N\to\infty} \Prob\biggl(
\frac{ \ell_N - 2\sqrt{N}}{N^{1/6}} \le x\biggr)
= F(x),
\end{equation}
where again $F(x)$ is the Tracy-Widom distribution appearing in 
\eqref{e-CLT}.
The authors in \cite{BDJ} also proved the convergence of moments: 
if $\chi_N(\sigma) := \frac{ \ell_N - 2\sqrt{N}}{N^{1/6}}$ 
and $\chi$ is a random variable with distribution function $F(x)$, then 
for $m=0,1,2,\cdots$, 
\begin{equation}
  \lim_{N\to\infty} \Exp( \chi_N^m) = \Exp (\chi^m).
\end{equation}
In earlier work other authors proved large deviation results for $\ell_N$, 
\begin{eqnarray}
\label{e-ld3}
  \lim_{N\to\infty} \frac1{N} \log  \Prob\bigl(
\ell_N \le \sqrt{N}(2-y) \bigr) = -H(y), \\
\label{e-ld4}
  \lim_{N\to\infty} \frac1{\sqrt{N}} \log  \Prob\bigl(
\ell_N \ge \sqrt{N}(2+y) \bigr) = -I(y),
\end{eqnarray}
for $y>0$, where $H(y), I(y)$ are 
certain explicit positive functions.
The result \eqref{e-ld3} is due to Deuschel and Zeitouni \cite{DeZe2}
and the result \eqref{e-ld4} is essentially due to 
Sepp\"al\"ainen \cite{Se2}.
 
In two recent papers, \cite{Merkl2} \cite{Merkl1}, the authors 
have considered $\ell_N$ in the moderate deviation regime.
More precisely, for $0<\alpha<\frac13$, they showed \cite{Merkl2} that 
for $y>0$
\begin{equation}\label{e-mld2}
  \lim_{N\to\infty}
\frac{\log  \Prob\bigl(
\ell_N \le (2-yN^{-\alpha})\sqrt{N} \bigr)}
{y^3N^{1-3\alpha}}  = -\frac1{12},
\end{equation}
and \cite{Merkl1}
\begin{equation}\label{e-mld1}
  \lim_{N\to\infty} 
\frac{\log  \Prob\bigl(
\ell_N \ge (2+yN^{-\alpha})\sqrt{N} \bigr)}
{y^{3/2}N^{(1-3\alpha)/2}}    = -\frac43.
\end{equation}
These moderate deviation results can be motivated by noting that 
\begin{equation}
  F(x) \sim e^{x^{3}/12}
\qquad \text{as $x\to -\infty$},
\end{equation}
and 
\begin{equation}
  1-F(x) \sim \frac{e^{-(4/3)x^{3/2}}}{16\pi x^{3/2}}
\qquad \text{as $x\to +\infty$}.
\end{equation}
Thus from \eqref{e-BDJ}, one anticipates that as $N\to\infty$, 
\begin{equation}\label{e-ldif2}
\begin{split}
  \Prob\bigl( \ell_N \le (2-yN^{-\alpha})\sqrt{N} \bigr)
&= \Prob\bigl( \ell_N \le 2\sqrt{N}-(yN^{1/3-\alpha})N^{1/6} \bigr)  \\
&\sim \log F(-yN^{1/3-\alpha}) \sim -\frac1{12}y^3N^{1-3\alpha},
\end{split}
\end{equation}
and 
\begin{equation}\label{e-ldif1}
\begin{split}
  \Prob\bigl( \ell_N \ge (2+yN^{-\alpha})\sqrt{N} \bigr)
&= \Prob\bigl( \ell_N \ge 2\sqrt{N}+(yN^{1/3-\alpha})N^{1/6} \bigr)  \\
&\sim \log (1-F(yN^{1/3-\alpha})) \sim -\frac43 y^{3/2}N^{(1-3\alpha)/2}.
\end{split}
\end{equation}
Of course, when $\alpha=0$, we are in the large deviation regime, and when 
$\alpha=\frac13$, we are in the GUE central limit theorem regime.
The above moderate deviation results can also be motivated 
by estimating the functions $I(y)$ and $H(y)$ for the large-deviation 
regime.
In \cite{Merkl1}, the authors proved \eqref{e-mld1} by refining 
certain estimates in \cite{BDJ} 
and using a careful summation argument. 
In \cite{Merkl2}, the authors utilized 
an analogous summation argument
together with the estimate Lemma 6.3 (ii) in \cite{BDJ}.

Calculations similar to \eqref{e-ldif2}, \eqref{e-ldif1}, 
motivate the following moderate deviation 
results for $G([\gamma N], N)$: for $0<\alpha<\frac23$, 
\begin{equation}\label{e-mld4}
  \lim_{N\to\infty}
\frac{\log  \Prob\bigl(
G([\gamma N],N) \le (\frac1{a_0}-yb_0N^{-\alpha})N \bigr)}
{y^3N^{2-3\alpha}}  = -\frac1{12},
\end{equation}
and 
\begin{equation}\label{e-mld3}
  \lim_{N\to\infty}
\frac{\log  \Prob\bigl(
G([\gamma N],N) \ge (\frac1{a_0}+yb_0N^{-\alpha})N \bigr)}
{y^{3/2}N^{1-3\alpha/2}}    = -\frac43.
\end{equation}
One of the principal goals in this paper is to prove \eqref{e-mld4}.
Relation \eqref{e-mld3} is slightly simpler and 
can also be approached using the techniques in this paper.
We hope to return to this problem in the future.

Relation \eqref{e-mld4} is a consequence of the following result.

\begin{thm}\label{mainthm}
Fix $0<t<1$ and $\gamma_0\ge 1$.
Then there exist a (large) constant $L>0$ and a (small) constant $\delta>0$,
such that for large $N$,
\begin{equation}\label{e-main}
  \log \Prob\bigl( G([\gamma N], N) \le \frac1{a_0} N  - x b_0 N^{1/3}\bigr)
= -\frac1{12} x^3 + O(x^4N^{-2/3}) + O(\log x)
\end{equation}
uniformly for all $L\le x\le \delta N^{2/3}$ and $1\le \gamma \le \gamma_0$.
In particular, for the variables $x, \gamma$ in the same range,
\begin{equation}\label{e-mainbound}
  \Prob\bigl( G([\gamma N], N) \le \frac1{a_0} N  - x b_0 N^{1/3}\bigr)
\le e^{-c|x|^3}
\end{equation}
for some constant $c>0$.
\end{thm}

Setting $x=yN^{2/3-\alpha}$ in \eqref{e-main}, 
we immediately obtain 
\eqref{e-mld4} together with error estimates.

\begin{cor}[Estimate for lower moderate deviation]\label{maincor}
For $0<\alpha<\frac23$ and $y>0$, as $N\to\infty$,
\begin{equation}\label{e-main2}
\frac{\log  \Prob\bigl(
G([\gamma N],N) \le (\frac1{a_0}-yb_0N^{-\alpha})N \bigr)}
{y^3N^{2-3\alpha}}  = -\frac1{12}
+O(yN^{-\alpha})+O\bigl(\frac{\log(yN^{2/3-\alpha})}{y^3N^{2-3\alpha}}\bigr),
\end{equation}
\end{cor}

Theorem \ref{mainthm} can also be used for other applications.

\begin{cor}[Convergence of moments]\label{cor-moments}
For $\gamma\ge 1$, 
set $\theta_N:= \frac{G([\gamma N], N)-a_0^{-1}N}{b_0N^{1/3}}$ 
and let $\chi$ be the random variable with distribution function 
$F(x)$ as above.
Then for $m=0,1,2,\cdots$, 
\begin{equation}\label{e-moments}
  \lim_{N\to\infty} \Exp( \theta_N^m) = \Exp (\chi^m).
\end{equation}
\end{cor}

\begin{proof}
By Remark 2.5 of \cite{kurtj:shape}, \eqref{e-moments} follows from 
the estimate \eqref{e-mainbound}.
\end{proof}

In particular, setting $m=2$, we see that the fluctuation 
$\sqrt{Var(G([\gamma N], N))}$ of $G$ is of order 
$N^{\eta}$ where $\eta=\frac13$.
It is believed (see e.g., \cite{spohn}) 
that the transversal fluctuations of $G$ 
have order $N^\xi$ where $\xi$ and $\eta$ are related by a
dimension-independent universal scaling law $2\xi=\eta+1$. 
In other words, it is expected that in our case $\xi=\frac23$.
In \cite{kurtj:trans}, Johansson considered transversal fluctuations 
for the Poissonized version of the longest increasing subsequence problem
and showed in that case that the scaling law $2
\xi = \eta + 1$ is satisfied.
By \cite{BDJ}, $\eta$ is again $1/3$ and it follows therefore that
the scaling law
$2\xi=\eta+1$ is satisfied for this case.
A key role in his analysis was again played by Lemma 6.3 (ii) of \cite{BDJ}.
This Poissonized problem can be viewed as a continuum version of 
the above site percolation problem 
and in Remark 1.2 of \cite{kurtj:trans} Johansson notes that 
the scaling law $2\xi=\eta+1$ for the site percolation problem would follow  
from an estimate of type \eqref{e-mainbound} above.
The modifications in the argument in \cite{kurtj:trans} that 
are needed for the site percolation problem are detailed 
in \cite{kurtj:pers}.
We thus have 

\begin{cor}[Transversal fluctuations]\label{cor-trans}
For any $\gamma\ge 1$, the above coefficients 
$\eta$ and $\xi$ for longitudinal and 
transversal fluctuations of the site percolation model 
obey the scaling law 
\begin{equation}
  2\xi=\eta+1.
\end{equation}
\end{cor}

In order to prove Corollary \ref{cor-moments}, \ref{cor-trans}, 
weaker bounds than \eqref{e-mainbound} suffice.
Indeed, using an observation of Harold Widom \cite{Widommom}
(see in particular Lemma 2), it is possible to prove the 
bound
\begin{equation}\label{e-weakbound}
  \Prob\bigl( G([\gamma N], N) \le \frac1{a_0} N  - x b_0 N^{1/3}\bigr)
\le e^{-c'|x|^{3/2}}
\end{equation}
for $x, \gamma$ in the range of Theorem \ref{mainthm},
for some constant $c'>0$.
As opposed to the proof of \eqref{e-mainbound}, which requires
a steepest-descent Riemann-Hilbert analysis (see below), 
the proof of \eqref{e-weakbound} uses only classical steepest-descent 
methods.
A key role in \cite{Widommom} is played by a beautiful conjecture 
of Widom for the spectral properties of a class of singular 
integral operators (see \cite{Widommom}).
This conjecture can be verified in our case, as in the case 
considered considered by Widom in \cite{Widommom}, 
by using an elegant formula 
of Borodin and Okounkov \cite{BoOk} (see identity (4.9)).
The method in \cite{Widommom} is itself motivated by earlier 
calculations in \cite{BDR}.
The estimate \eqref{e-weakbound} is enough to prove 
Corollary \ref{cor-moments}, \ref{cor-trans}, but does not suffice to prove 
Corollary \ref{maincor}.

\begin{rem}
 The bound in Lemma 6.3 (ii) of \cite{BDJ} was also used 
by Sepp\"al\"ainen \cite{Seppal} to control fluctuations for
the ``stick process'' introduced in \cite{Se95}.
In \cite{Seppal}, Sepp\"al\"ainen also mentioned that 
a similar result could be obtained for a certain 
continuous-time totally asymmetric 
simple exclusion process, provided the appropriate analogue 
of Lemma 6.3 (ii) could be established.
The same should be true for a discrete-time version of this 
process.
The above estimate \eqref{e-weakbound}, and of course also the
stronger estimate \eqref{e-mainbound}, then suffices to control 
the fluctuations as in the stick process.
\end{rem}

\begin{rem}
   Our results are given for $G(M,N)$ where $M=[\gamma N]\ge N$, 
but as the statistics of $G(M,N)$ are clearly the same as for $G(N,M)$, 
it is immediate that our results, suitably scaled, also apply 
to $G([\gamma N], N)$ for $0<\gamma <1$.
\end{rem}

As indicated above, the proof of Theorem \ref{mainthm} 
is based on the steepest-descent method 
for Riemann-Hilbert problems (RHP's) introduced by Deift and Zhou 
\cite{DZ1}, and further developed in \cite{DVZ}.
The method has been used to solve a wide variety of asymptotic problems 
in pure and applied mathematics (see, for example, \cite{DKMVZ2, DKMVZ3} 
and the references therein).
The steepest-descent calculations in this paper are closely related to 
the calculations in \cite{DKMVZ2, DKMVZ3} and 
particularly \cite{BDJ}.
Our analysis is based on 
the algebraic formula \eqref{e-Gprod} below, 
which relates $\Prob(G(M,N) \le n)$ to the solution $Y$ to an 
associated RHP on the unit circle $\Sigma=\{|z|=1\}$ (see \eqref{e-Y}).
It follows then that our problem reduces to the asymptotic analysis 
of a RHP with large oscillatory parameters.
The steepest-descent method in \cite{DZ1} 
was introduced precisely for this purpose.
The key step in the method is to identify the leading order asymptotics 
for the solution of the RHP and this is done following 
\cite{DVZ}, \cite{DKMVZ2, DKMVZ3} and \cite{BDJ} 
by introducing a so-called $g$-function with certain specific properties 
on an appropriate contour $\overline{\Gamma_1\cup\Gamma_2}$
(see Proposition \ref{prop-g} below).
Using $g$, one transforms the RHP for $Y$ as follows: 
$U:= e^{-\frac12k\ell\sigma_3}Ye^{-k(g-\frac12\ell)\sigma_3}$ 
where $\ell$ is a specific constant to be determined and 
$\sigma_3$ is the Pauli matrix $\sigma_3= \bigl( \begin{smallmatrix}
1&0\\0&-1 \end{smallmatrix} \bigr)$.
A simple calculation shows that $U$ solves the RHP 
\begin{equation}\label{e-U}
\begin{cases}
  U(z) \text{ is analytic in $z\in\C\setminus \Sigma$,
and continous up to the boundary}, \\
  U_+(z)= U_-(z) \begin{pmatrix} e^{-k(g_+-g_-)} 
& e^{k(g_++g_- -W-\ell)} \\ 0 & e^{k(g_+-g_-)} \end{pmatrix},
\qquad z\in\Sigma, 
\end{cases}
\end{equation}
where $W$ is given \eqref{e-Yagamma}, 
and  $g_\pm$ denote the boundary values of $g$.
In addition, one requires $g(z)= \log z +O(z^{-1})$ as $z\to\infty$, 
so that the RHP for $U$ is normalized at infinity,
\begin{equation}
  U(z) = I + O(1/z), \qquad \text{as $z\to\infty$}.
\end{equation}
The choice of the properties of $g$ mentioned above is made precisely 
such that the leading contribution to the RHP \eqref{e-U} 
is immediate. Further information on the steepest-descent method 
can be found, for example, in \cite{DVZ}, \cite{DKMVZ2, DKMVZ3}, \cite{BDJ}.

In \cite{DVZ}, \cite{DKMVZ2, DKMVZ3}, the RHP's are given 
on the real line $\R$ and the analogues of $\Gamma_1$, $\Gamma_2$ 
are subintervals of $\R$: in \cite{BDJ}, the RHP 
is given on the unit circle $\Sigma=\{|z|=1\}$ and the analogues 
of $\Gamma_1$, $\Gamma_2$
are again subintervals of $\Sigma$.
The main new technical feature of the RHP in this paper 
is that $\Gamma_1$ and $\Gamma_2$ cannot be chosen as subintervals 
of the original contour $\Sigma$, and 
the central problem is to discover the shape of $\Gamma_1$, $\Gamma_2$.
The situation is similar to the problem confronted by 
Kamvissis, McLaughlin and Miller in \cite{KMM}, where the authors 
considered the semi-classical limit of the solution of the Cauchy 
problem for the focusing nonlinear Schr\"odinger equation.
Motivated by the calculations in \cite{KMM}, we show that the construction 
of $\Gamma_1$, $\Gamma_2$ is equivalent to the problem of determining 
the global structure of the trajectories $Q(z)(dz)^2>0$ 
and orthogonal trajectories $Q(z)(dz)^2<0$ of a particular 
quadratic differential $Q(z)(dz)^2$ (see \eqref{e-Q} below).

The outline of the paper is as follows.
In Section \ref{sec-setup}, we derive 
the basic algebraic formula \eqref{e-Gprod} 
relating $\Prob(G([\gamma N], N) \le n)$ and the RHP \eqref{e-Y}, 
and state our basic asymptotic estimate, 
Proposition \ref{prop1}, for $Y_{21}(0;k)$.
In Section \ref{sec-h}, which is the heart of the paper,  
we construct $\Sigma_1$ and $\Sigma_2$ using the theory of 
quadratic differentials 
and verify the desired properties of $h=g'$.
In Section \ref{sec-g},
the constant $\ell$ mentioned above is defined (see \eqref{e-l}) and 
the desired properties of $g=\int^z h$ are verified 
(Proposition \ref{prop-g}).
In Section \ref{sec-rhp}, we use the $g$-function to analyze 
the RHP and eventually give the proof of Proposition \ref{prop1}.
Finally, in Section \ref{sec-proof}, we use the estimate in  
Proposition \ref{prop1} together with a careful summation argument 
as in \cite{Merkl2} to prove the main result Theorem \ref{mainthm}.

\medskip
\noindent {\bf Acknowledgments.} The authors would like to thank
Nick Ercolani for useful discussions and Kurt Johansson 
for making available to us his calculations in \cite{kurtj:pers}.
The authors would also like to thank Harold Widom for providing 
us with his preprint \cite{Widommom}.
The first author would like to thank Anne Boutet de Monvel for 
kindly inviting him to Universit\'e Paris 7, where a part 
of work is done, and also acknowledge that 
a part of work is conducted while he is visiting 
Korea Institute for Advanced Study for 2 weeks of August, 2001.
The work of the first author was
supported in part by NSF Grant \# DMS 97-29992. The work of the
second author was supported in part by NSF Grant \# DMS 00-03268.
The work of the third author was supported in part by NSF Grant \#
DMS-9970328.
The work of the fourth author was supported in part by NSF Grant \#
DMS 01-03909.
The work of the fifth author was supported in part by NSF Grant \#
DMS 0071398.

\section{Basic relations and formulae}\label{sec-setup}

For $M,N \ge 1$, let 
\begin{equation}
  Z_{M,N}:= (1-t^2)^{-MN}.
\end{equation}
Set 
\begin{equation}
  \varphi(z)= (1+tz)^M(1+\frac{t}{z})^N,
\end{equation}
and consider the $n\times n$ Toeplitz determinant
\begin{equation}
 D_n(\varphi)=D_n := \det( \varphi_{j-k})_{0\le j,k <n},
\end{equation}
where $\varphi_j$ is the $j^{th}$ Fourier coefficient of $\varphi$: 
\begin{equation}
  \varphi_j:= \int_{|z|=1} z^{-j} \varphi(z) \frac{dz}{2\pi iz}.
\end{equation}
Here and below the integration contour $|z|=1$ is assumed to be
oriented in the counter-clockwise direction.  Let $G(M,N)$ be the
maximal length introduced in the Introduction.  From earlier result
of Gessel \cite{Gessel} an Johansson \cite{kurtj:shape}, Baik and
Rains \cite{BR1} extracted the relation 
\begin{equation}\label{e-Gdet}
  \Prob( G(M,N) \le n)  = \frac1{Z_{M,N}} D_n(\varphi),
\end{equation}
which plays the basic role in our analysis.

Let $\Sigma$ be the unit circle $|z|=1$ in the complex plane, 
oriented counter-clockwise
and let $Y(z)=Y(z;k)= (Y_{ij}(z;k))_{1\le i,j \le2}$ 
be the solution to the following $2\times 2$ matrix RHP: 
\begin{equation}\label{e-Y}
\begin{cases}
  Y(z) \text{ is analytic in $z\in\C\setminus \Sigma$, 
and continous up to the boundary}, \\
  Y_+(z)= Y_-(z) \begin{pmatrix} 1 & z^{-k}\varphi(z) \\ 0 & 1 \end{pmatrix},
\qquad z\in\Sigma, \\
  Y(z)z^{-k\sigma_3} = I + O(1/z), \qquad \text{as $z\to\infty$},
\end{cases}
\end{equation}
where $\sigma_3= \bigl( \begin{smallmatrix} 1&0 \\ 0&-1 \end{smallmatrix}
\bigr)$ is the standard third Pauli matrix, 
and $Y_+(z)$, (resp., $Y_-(z)$), 
$z\in\Sigma$, are the boundary values of $Y(z')$ as $z'\to z$
from the inside (resp., outside) of the circle.

\begin{lem}\label{lem-GandY}
The solution $Y$ to the above RHP \eqref{e-Y} exists and is unique. Moreover, 
\begin{equation}\label{e-Gprod}
  \Prob( G(M,N) \le n) = \prod_{k=n}^{\infty} (-Y_{21}(0; k+1)).
\end{equation}
\end{lem}

\begin{proof}
We will construct the solution $Y$ explicitly 
using computations similar to \cite{FIK}.
First note that from the equality \eqref{e-Gdet}, $D_n(\varphi)\neq 0$ for 
$n\ge 0$ since the probability 
$\Prob(G(M,N)=0)= \Prob(w(i,j)=0,  
1\le \forall i\le M, 1\le \forall j\le N) = (1-q)^{MN}$, and hence 
$\Prob(G(M,N)\le n)\ge \Prob(G(M,N)=0)>0$ for $n\ge 0$.
Consider for $k\ge 0$ the polynomials of degree $k$
\begin{equation}\label{e-1.8}
  \pi_k(z):= \frac1{D_{k}} \det \begin{pmatrix}
\varphi_0 & \varphi_{-1} & \cdots &\varphi_{-k} \\
\varphi_{1} & \varphi_0 & \cdots &\varphi_{-k+1} \\
\vdots & \vdots & & \vdots \\
\varphi_{k-1} & \varphi_{k-2} &\cdots & \varphi_{-1} \\
1 & z & \cdots & z^k
\end{pmatrix},
\quad
  \pi^*_k(z):= \frac1{D_{k}} \det \begin{pmatrix}
\varphi_0 & \varphi_{1} & \cdots &\varphi_{k} \\
\varphi_{-1} & \varphi_0 & \cdots &\varphi_{k-1} \\
\vdots & \vdots & & \vdots \\
\varphi_{-k+1} & \varphi_{-k+2} &\cdots & \varphi_{1} \\
z^k & z^{k-1} & \cdots & 1
\end{pmatrix}.
\end{equation}
A direct check shows that $\pi_k$ and $\pi^*_k$ satisfy 
the following orthogonality conditions:
\begin{eqnarray}
\label{e-orth1}
  \int_{|z|=1} z^{-j} \pi_k(z) \varphi(z) \frac{dz}{2\pi iz} 
&=& N_k \delta_{jk}
\qquad 0\le j\le k, \\
\label{e-orth2}
  \int_{|z|=1} z^{-j} \pi^*_k(z) \varphi(z) \frac{dz}{2\pi iz} 
&=& N_k \delta_{j0}
\qquad 0\le j\le k.
\end{eqnarray}
where 
\begin{equation}\label{e-1.11}
  N_k= \frac{D_{k+1}}{D_k}.
\end{equation}
Let $(Ch)(z)= \frac1{2\pi i}\int_{\Sigma} \frac{h(s)}{s-z} ds$, 
$z\in \C\setminus\Sigma$, denote the Cauchy transform of $h$.
Let $(C_\pm h)(z)= \lim_{z'\to z} (Ch)(z')$ 
where $z'$ approaches $z$ from the $\pm$ side respectively, 
denote its boundary values as in the Introduction.
We claim that 
\begin{equation}\label{e-Ysol}
  Y(z; k) = \begin{pmatrix}
\pi_k(z) & C((\cdot)^{-k}\varphi\pi_k)(z) \\ 
- N_{k-1}^{-1} \pi^*_{k-1}(z) & 
-N_{k-1}^{-1} C((\cdot)^{-k}\varphi\pi^*_{k-1})(z)
\end{pmatrix}.
\end{equation}
is a solution to \eqref{e-Y}.
The analyticity 
of $Y$ in $\C\setminus\Sigma$ is clear, while the jump 
condition follows from the relation $C_+-C_-=1$.
The asymptotic condition follows from the orthogonality 
\eqref{e-orth1}, \eqref{e-orth2}.
On the other hand, the uniqueness of the solution to the RHP \eqref{e-Y} 
is standard (cf. for example, Lemma 4.1 of \cite{BDJ}).
Hence \eqref{e-Ysol} is the unique solution to the RHP \eqref{e-Y}.

For the proof of \eqref{e-Gprod}, first note that by taking $n\to\infty$ in 
\eqref{e-Gdet}, 
\begin{equation}
  \lim_{n\to\infty} D_n(\varphi) = Z_{M,N}.
\end{equation}
(This can also be seen directly from 
the Szeg\"o strong limit theorem for Toeplitz 
determinants.)
Thus we have, using \eqref{e-1.11},  
\begin{equation}
 \Prob( G(M,N) \le n) = \prod_{k=n}^\infty \frac{D_k}{D_{k+1}} 
= \prod_{k=n}^\infty N_k^{-1}
\end{equation}
Finally, from \eqref{e-1.8} and \eqref{e-Ysol}, we observe that 
$Y_{21}(0; k)= - N_{k-1}^{-1}$, which completes the proof.
\end{proof}

\begin{rem}
From \eqref{e-Gprod}, 
\begin{equation}
  -Y_{21}(0;k+1) = \frac{\Prob(G(M,N)\le k)}{\Prob(G(M,N)\le k+1)},
\end{equation}
and hence we have $-Y_{21}(0;k) >0$ for $k\ge 1$.
\end{rem}

\bigskip

There are three parameters in the RHP \eqref{e-Y}: $M, N, k$.
We regard $t$, $0<t<1$, as a fixed number throughout this paper.
For convenience we introduce the following notation:
\begin{equation}
  \gamma := \frac{M}{N}, \qquad a:= \frac{N}{k}.
\end{equation}
Instead of $M, N, k$, we now regard $\gamma, a, k$ as our parameters 
in the RHP \eqref{e-Y}.
Using this notation, the jump condition 
for $Y$ takes the form 
\begin{equation}\label{e-Yagamma}
  Y_+(z)= Y_-(z) \begin{pmatrix} 
1& e^{- k W(z)}
\\ 0&1 \end{pmatrix}, 
\qquad W(z):= - \gamma a\log(1+tz) - a\log(1+t/z) + \log z,
\end{equation}
where $\log w$ is defined to be analytic in $\C\setminus (-\infty, 0]$, 
and $\log w = \log|w|$ for $w>0$.  
We are interested in the asymptotics of $Y$ as $k\to\infty$, while 
$\gamma$ and $a$ remain in appropriate bounded regions.
The main technical result we are going to prove 
in the rest of this paper is the following:

\begin{prop}\label{prop1}
Set 
\begin{equation}\label{e-a0}
  a_0= \frac{1-t^2}{t((\gamma+1)t + 2\sqrt{\gamma})}.
\end{equation}
Fix $0<t<1$ and $\gamma_0\ge 1$.
There are $L_0$, $\delta_0$, $k_0 >0$ such that 
for $\gamma$, $a$, $k$ satisfying 
\begin{equation}\label{e-kcond}
  a_0 + \frac{L_0}{k^{2/3}} \le a \le (1+\delta_0) a_0,
\end{equation}
$k\ge k_0$ and $1\le \gamma \le \gamma_0$, we have 
%
%
\begin{equation}\label{e-Yexactest}
 \log (-Y_{21}(0; k)) = -c_2k(a-a_0)^2 + O(k|a-a_0|^3)+O(|a-a_0|)
+O\bigl(\frac1{k|a-a_0|}\bigr),
\end{equation}
where 
\begin{equation}
  c_2=\frac{t^2(t+t\gamma+2\sqrt\gamma)^3\sqrt\gamma}
{4(1+t\sqrt\gamma)^2(t+\sqrt\gamma)^2}.
\end{equation}
In particular, for the variables in the same range, 
\begin{equation}\label{e-Yineq}
  \log (-Y_{21}(0; k)) \le -c_0 k(a-a_0)^2
\end{equation}
for some constant $c_0>0$.
\end{prop}

\section{The $h$-function}\label{sec-h}
\begin{figure}[ht]
 \centerline{\epsfig{file=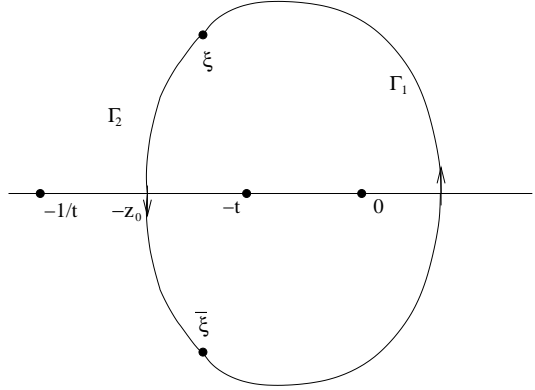, width=7cm}}
 \caption{Sketch of the contours $\Gamma_1$ and $\Gamma_2$.}
\label{fig-fig0}
\end{figure}

As noted in the Introduction,  we seek a $g$-function and
associated contours  $\overline{\Gamma_1\cup\Gamma_2}$. Rather
than analyzing $g$ directly, we first seek  an \emph{$h$-function},
$h'=g$. We start from an ansatz, to be verified a posteriori, that
$\Gamma_1$, $\Gamma_2$  has the shape  indicated in Figure
\ref{fig-fig0}  with the endpoints $\xi$, $\overline{\xi}$.
Recall from \eqref{e-Yagamma}
\begin{equation}\label{e-W}
  W(z) = -\gamma a\log(1+tz) - a\log(1+t/z) + \log z, 
\qquad W'(z)= -\frac{\gamma a}{z+t^{-1}}  - \frac{a}{z+t} + \frac{a+1}{z},
\end{equation}
and recall from \eqref{e-a0}
\begin{equation}
  a_0 = \frac{1-t^2}{t((\gamma+1)t+2\sqrt\gamma)}.
\end{equation}

\noindent {\it Notation from the theory of quadratic differentials 
(see e.g., \cite{Pomm}):}
Given a meromorphic function $f(z)$ and 
a simple oriented curve $C$ lying outside
the zeros and poles of $f$, 
the notation $f(z)(dz)^2>0$ means that 
$f(z(t))(\frac{dz}{dt})^2$ is real and positive for all $t\in (a,b)$ 
where $z(t)$, $t\in (a,b)$ is the arc length parameterization 
of $C$.
Similarly, $f(z)(dz)^2<0$ means that 
$f(z(t))(\frac{dz}{dt})^2$ is real and negative.

\begin{prop}\label{prop-h}
Fix $0<t<1$, $\gamma\ge 1$, $a>a_0$.
Then there exist (cf. Figure \ref{fig-fig0} below) 
\begin{itemize}
\item a point $\xi$ with $Im(\xi)>0$,
\item a simple open curve $\Gamma_1$ connecting
$\overline{\xi}$ and $\xi$, oriented from $\overline{\xi}$ to $\xi$ such that
(i) it does not intersect $(-\infty, 0]$
and (ii) $\Gamma_1$ is symmetric with reflection about the real line,
\item a simple open curve $\Gamma_2$ connecting
$\xi$ and $\overline{\xi}$, oriented from $\xi$ to $\overline{\xi}$
such that (i) it does not intersect $(-\infty, -1/t]\cup [-t,\infty)$
and (ii) $\Gamma_2$ is symmetric with respect to reflection about the
real line,
\item a function $h(z)$ analytic in
$\C\setminus\overline{\Gamma_1}$ and
continuous up to the boundary
\end{itemize}
such that the following properties are satisfied:
\begin{itemize}
\item[(a)] $h_+(z)+h_-(z) = W'(z)$ for $z\in\Gamma_1$.  
\item[(b)] $h(z) = \frac1{z}+ O(z^{-2}))$ as $z\to\infty$.
\item[(c)] $i(h_+(z)-h_-(z))dz >0$ for $z\in\Gamma_1$.
\item[(d)] $(2h(z)- W'(z))dz <0$ for $z\in\Gamma_2\cap\C_+$ 
and $(2h(z)- W'(z))dz >0$ for $z\in\Gamma_2\cap\C_-$.
\end{itemize}
In addition, we have the following properties: 
\begin{itemize}
\item[(i)] The function $h$ has the form 
\begin{equation}\label{e-hprop}
\begin{split}
  h(z) &= \frac{R(z)}{2\pi i} \int_{\Gamma_1} \frac{W'(s)}{R_+(s)(s-z)} ds \\
&= \frac12 W'(z) + \frac12 R(z)\biggl(
\frac{\gamma a}{(z+t^{-1})R(-t^{-1})} + \frac{a}{(z+t)R(-t)}
-\frac{a+1}{zR(0)} \biggr), 
\end{split}
\end{equation}
where $R(z)=\sqrt{(z-\xi)(z-\overline{\xi})}$ is defined to be 
analytic in $\C\setminus\Gamma_1$ and $R(z)\sim z$ as $z\to\infty$.
\item[(ii)] (Endpoint condition) 
\begin{equation}\label{e-endprop}
  \int_{\Gamma_1} \frac{W'(s)}{R_+(s)} ds =0, \qquad
\frac1{2\pi i} \int_{\Gamma_1} \frac{sW'(s)}{R_+(s)} ds =-1.
\end{equation}
\item[(iii)]
Set 
\begin{equation}\label{e-Phi}
  \Phi(z) := (1+\gamma a) \frac{R(z)(z+z_0)}{z(z+t)(z+t^{-1})}, 
\qquad z\in\C\setminus\Gamma_1,
\end{equation}
where 
\begin{equation}\label{e-z0inprop}
  z_0= \frac{a+1}{-R(0)(1+\gamma a)}.
\end{equation}
Then 
\begin{equation}\label{e-hPhi}
  2h(z)-W'(z)= \Phi(z), \qquad z\in\C\setminus\overline{\Gamma_1}.
\end{equation}
In particular, 
\begin{eqnarray}
\label{e-hPhi1}
  &h_+(z)-h_-(z) = \Phi_+(z),&  \qquad z\in\Gamma_1, \\
\label{e-hPhi2} 
  &2h(z)-W'(z) = \Phi(z),& \qquad z\in\Gamma_2.
\end{eqnarray}
\item[(iv)] The curve $\Gamma_2$ intersects the real axis at $-z_0$.
\end{itemize}
\end{prop}

\noindent{\it Idea of proof:}

Suppose that the curve $\Gamma_1$ is known.
Let
\begin{equation}\label{e-R-0}
R(z)=\sqrt{(z-\xi)(z-\overline{\xi})},
\end{equation}
which is defined to be
analytic in $\C\setminus\Gamma_1$, and $R(z)\sim z$ as $z\to\infty$.
From the Plemelj formula, it is easy to check that
\begin{equation}\label{e-h}
  h(z) = \frac{R(z)}{2\pi i} \int_{\Gamma_1} \frac{W'(s)}{R_+(s)(s-z)} ds,
\end{equation}
satisfies condition (a).
By a residue calculation, $h$ can be written as
\begin{equation}\label{e-hres-0}
  h(z) = \frac12 W'(z) + \frac12 R(z)\biggl(
\frac{\gamma a}{(z+t^{-1})R(-t^{-1})} + \frac{a}{(z+t)R(-t)}
-\frac{a+1}{zR(0)} \biggr).
\end{equation}
For this $h$ to satisfy (b), 
the following two conditions are necessary and sufficient,
\begin{equation}\label{e-end-0}
  \int_{\Gamma_1} \frac{W'(s)}{R_+(s)} ds =0, \qquad
\frac1{2\pi i} \int_{\Gamma_1} \frac{sW'(s)}{R_+(s)} ds =-1,
\end{equation}
or equivalently by residue calculations, 
\begin{equation}\label{e-endcon-0}
  \frac{\gamma a}{R(-t^{-1})} + \frac{a}{R(-t)} = \frac{a+1}{R(0)},
\qquad  \frac{\gamma a}{tR(-t^{-1})} + \frac{at}{R(-t)} = -(1+\gamma a).
\end{equation}
As we will see, these conditions determine the endpoint $\xi$.

Now for $z\in\Gamma_1$, from \eqref{e-hres-0}, 
\begin{equation}\label{e-hdiff-0}
  i(h_+(z)-h_-(z)) = iR_+(z) \biggl(
\frac{\gamma a}{(z+t^{-1})R(-t^{-1})} + \frac{a}{(z+t)R(-t)}
-\frac{a+1}{zR(0)} \biggr).
\end{equation}
Substituting $R(-t^{-1})$ and $R(-t)$ in terms of $R(0)$ using the
endpoint conditions \eqref{e-endcon-0}
we obtain, after some algebra,  for $z\in\Gamma_1$,
\begin{equation}
  i(h_+(z)- h_-(z)) = i(1+\gamma a) R_+(z) \frac{(z+z_0)}{z(z+t)(z+t^{-1})},
\end{equation}
where
\begin{equation}
  z_0 = \frac{a+1}{-R(0)(1+\gamma a)}.
\end{equation}
Hence we have 
\begin{equation}\label{e-Qdef1}
  Q(z):=(i(h_+-h_-))^2 = -(1+\gamma a)^2
\frac{(z-\xi)(z-\overline{\xi})(z+z_0)^2}{z^2(z+t)^2(z+t^{-1})^2},
\end{equation} 
which is meromorphic in $\C$.
We then use the theory of the quadratic differentials 
to find the trajectories for $Q(z)(dz)^2>0$ and this leads us to the 
determination of the contour $\Gamma_1$ for which 
condition (c) is satisfied.

The contour $\Gamma_2$ for which condition (d) is satisfied 
turns out to be obtained by finding the 
so-called orthogonal trajectories corresponding to $Q(z)(dz)^2<0$.
Clearly if a trajectory and an orthogonal trajectory meet
at a point $z$ in the plane where $Q(z)$ is analytic and nonzero, they do 
so at right angles to each other.

\bigskip

The rest of this section consists of a proof 
of the above Proposition.

\subsection{The endpoint $\xi$ }\label{subsec-end}

In this subsection, 
we are going to prove that there is a unique $\xi$, $Im(\xi)>0$, 
for which the following two conditions (cf. \eqref{e-end-0})
are satisfied,
\begin{equation}\label{e-end1}
  \int_{\Gamma} \frac{W'(s)}{R_+(s)} ds =0, \qquad 
\frac1{2\pi i} \int_{\Gamma} \frac{sW'(s)}{R_+(s)} ds =-1,
\end{equation}
where $\Gamma$ is \emph{any} simple oriented curve connecting
$\xi$ and $\overline{\xi}$, oriented from $\overline{\xi}$ to $\xi$
such that (i) it does not intersect $(-\infty, 0]$ 
(ii)  it is symmetric under reflection about the real line, 
and 
\begin{equation}\label{e-R}
R(z)=\sqrt{(z-\xi)(z-\overline{\xi})},
\end{equation}
which is defined to be
analytic in $\C\setminus\Gamma$ with $R(z)\sim z$ as $z\to\infty$.

\begin{rem}
A priori we should look for a pair of unrelated points $\xi_1$, $\xi_2$
such that \eqref{e-end1} is satisfied for any contour $\Gamma$ connecting
them as above.
It turns out, however, that it is sufficient to look for 
$\xi_1$ and $\xi_2$ in the form 
$\xi_1=\xi$ and $\xi_2=\overline{\xi}$.
The reason for this symmetry lies in the form of the equations
\eqref{e-end1}.  Indeed, both $W'(s)$ and $sW'(s)$ are real analytic
functions, and in each integral the path of integration can be doubled
along the ``minus'' side of the branch cut for $R(s)$ and then
deformed into a closed loop containing $\xi_1$ and $\xi_2$ and the
branch cut $\Gamma$ connecting them but no singularities of $W'(s)$;
this loop is otherwise arbitrary.  If we take the loop to be symmetric
with respect to reflection in the real axis, then it is easy to see
that the only way for both integrals to be purely imaginary as
required by \eqref{e-end1} is for $R(s)$ itself to be a real analytic
function, which forces $\xi_2=\overline{\xi_1}$.
\end{rem}

As in \eqref{e-endcon-0}, these conditions become 
\begin{equation}\label{e-endcon}
  \frac{\gamma a}{R(-t^{-1})} + \frac{a}{R(-t)} = \frac{a+1}{R(0)}, 
\qquad  \frac{\gamma a}{tR(-t^{-1})} + \frac{at}{R(-t)} = -(1+\gamma a).
\end{equation} 
Set 
\begin{equation}\label{e-rxyxi}
  r:= |\xi| =-R(0), \qquad x:=|t^{-1}+\xi|=-R(-t^{-1}), 
\qquad y:= |t+\xi|= - R(-t).
\end{equation}
Then conditions \eqref{e-endcon} have the form
\begin{equation}\label{e-contemp}
  \frac{\gamma a}{x} + \frac{a}{y} = \frac{a+1}{r}, 
\qquad \frac{\gamma a}{tx} + \frac{at}{y} = 1+\gamma a.
\end{equation}
Set 
\begin{equation}
  r_1:= \frac{t(a+1)}{1+\gamma a}, \qquad 
r_2 := \frac{a+1}{t(1+\gamma a)}.
\end{equation}
The conditions \eqref{e-contemp}
are now, after simple algebra, equivalent to 
\begin{eqnarray}
\label{e-con1}
  \frac{a\gamma (1-t^2)}{x} &=& t(1+\gamma a)\bigl(1 - \frac{r_1}{r}\bigr), \\
\label{e-con2}
\frac{a(1-t^2)}{y}&=& t(1+\gamma a) \bigl( \frac{r_2}{r}-1 \bigr).
\end{eqnarray}
Also from the definitions, $r, x, y$ satisfy the relation
\begin{equation}\label{e-con3}
  r^2 = 1+ \frac{y^2-t^2x^2}{1-t^2}.
\end{equation}
Inserting $x, y$ of \eqref{e-con1}, \eqref{e-con2} into 
\eqref{e-con3}, we obtain an equation for $r$: 
\begin{equation}\label{e-Heq}
  \frac1{r^2}  + \frac{(1-t^2)a^2}{(1+\gamma a)^2} 
\biggl( \frac1{t^2(r_2-r)^2} - \frac{\gamma ^2}{(r-r_1)^2} \biggr) -1 
=0.
\end{equation}
Since $x, y>0$, from \eqref{e-con1}, \eqref{e-con2}, we must have 
$r_1< r< r_2$. Thus we seek $x,y,r$ satisfying 
\begin{equation}\label{e-rxyineq}
  x>0, \qquad y>0,\qquad  r_1<r<r_2.
\end{equation}

\begin{lem}\label{lem-H1}
For each fixed $0<t<1$, $\gamma\ge 1$, $a>a_0$, 
there is a unique solution $r$ to \eqref{e-Heq} satisfying
$r_1< r <r_2$.
\end{lem}

\begin{proof}
Set 
\begin{equation}\label{e-H}
  H(r) := \frac1{r^2}  + \frac{(1-t^2)a^2}{(1+\gamma a)^2}
\biggl( \frac1{t^2(r_2-r)^2} - \frac{\gamma ^2}{(r-r_1)^2} \biggr) -1.
\end{equation}
Clearly, $H(r) \to -\infty$ as $r\downarrow r_1$, and 
$H(r) \to +\infty$ as $r\uparrow r_2$.
Thus there is $r_1<r_c<r_2$ satisfying $H(r_c)=0$.
We want to show that such an $r_c$ is unique.
By direct calculation, for $r_1<r<r_2$, 
\begin{equation}
  H(r) + \frac{r}2 H'(r) 
= \frac{a^2(1+a)(1-t^2)}{(1+\gamma a)^3} \biggl( 
\frac{1}{t^3(r_2-r)^3} + \frac{t\gamma^2}{(r-r_1)^3} \biggr)  -1.
\end{equation}
The minimum of this function on $(r_1, r_2)$ is obtained at
\begin{equation}\label{e-minr0}
  r_* = \frac{r_1+r_2 \sqrt\gamma t}{1+\sqrt\gamma t}, 
\end{equation}
and for $r_1< r< r_2$, 
\begin{equation}
  H(r)+ \frac{r}{2}H'(r) \ge H(r_*) + \frac{r_*}2H'(r_*)
= \biggl( \frac{a(1+\sqrt\gamma t)^2}{(1+a)(1-t^2)} \biggr)^2 -1.
\end{equation}
But since $a>a_0$, 
\begin{equation}
   H(r)+ \frac{r}{2}H'(r) >
\biggl( \frac{a_0(1+\sqrt\gamma t)^2}{(1+a_0)(1-t^2)} \biggr)^2 -1 
=0.
\end{equation}
Therefore, if $H(r_c)=0$ for $r_1<r_c<r_2$, we must have 
$H'(r_c)>0$. 
A simple calculus argument then proves the uniqueness of the solution.
\end{proof}

Thus if we define $x, y$ by \eqref{e-con1}, \eqref{e-con2}, 
we have obtained the unique solution $(r, x, y)$ 
to the equations \eqref{e-con1}, \eqref{e-con2}, \eqref{e-con3} 
subject to \eqref{e-rxyineq}.
Now we need to prove that the $(r, x, y)$ defined in this way determines 
$\xi$, $Im(\xi)>0$, uniquely from \eqref{e-rxyxi}.
In order for $\xi$, $Im(\xi)>0$, to satisfy \eqref{e-rxyxi}, 
we must have $\xi=re^{i\theta}$ for some 
$0<\theta<\pi$ satisfying 
\begin{equation}\label{e-cosxi}
  \cos\theta = \frac{x^2-r^2-t^{-2}}{2rt^{-1}} =
\frac{y^2-r^2-t^2}{2rt}.
\end{equation}
Conversely, if there exists $\theta \in (0,\pi)$ satisfying 
\eqref{e-cosxi}, then $\xi:=re^{i\theta}$ 
and $\overline{\xi}=re^{-i\theta}$ are the desired endpoints.
However, the second inequality follows from \eqref{e-con3} 
and so it is sufficient to prove that for $(x,y,r)$ 
satisfying \eqref{e-con1}, \eqref{e-con2}, \eqref{e-con3}, \eqref{e-rxyineq}, 
we have the relation 
\begin{equation}\label{e-xiextra}
  -1< \frac{x^2-r^2-t^{-2}}{2rt^{-1}} <1.
\end{equation}
In order to prove \eqref{e-xiextra}, we first prove the following Lemma.

\begin{lem}\label{lem-endptxyt}
For each $0<t<1$, $\gamma\ge 1$, $a>a_0$, 
the solution $(r, x, y)$ to \eqref{e-con1}, \eqref{e-con2}, 
\eqref{e-con3}, subject to \eqref{e-rxyineq} satisfies 
\begin{equation}
  x+y > t^{-1}-t.
\end{equation}
\end{lem}

\begin{proof}
From \eqref{e-con1}, \eqref{e-con2}, we have 
\begin{equation}
  \frac{x+y}{t^{-1}-t} = \frac{a}{1+\gamma a} \biggl( 
\gamma -1 + \frac{r_1\gamma}{r-r_1} + \frac{r_2}{r_2-r} \biggr).
\end{equation}
The minimum of the right-hand side, regarded as a function in $r$, 
is again obtained at $r=r_*$, 
where $r_*$ is defined in \eqref{e-minr0}, and hence, 
by evaluating the minimum, we obtain 
\begin{equation}
  \frac{x+y}{t^{-1}-t} \ge \frac{a(t+\sqrt\gamma)^2}{(1+\gamma a)(1-t^2)}. 
\end{equation}
But since $a>a_0$, we have 
\begin{equation}
  \frac{x+y}{t^{-1}-t} > \frac{a_0(t+\sqrt\gamma)^2}{(1+\gamma a_0)(1-t^2)}
=1.
\end{equation}
\end{proof}

Now we prove \eqref{e-xiextra}.

\begin{lem}
For each $0<t<1$, $\gamma\ge 1$, $a>a_0$, 
the solution $(r, x, y)$ to \eqref{e-con1}, \eqref{e-con2},
\eqref{e-con3}, subject to \eqref{e-rxyineq} satisfies 
\begin{equation}
  -1< \frac{x^2-r^2-t^{-2}}{2rt^{-1}} <1.
\end{equation}
\end{lem}

\begin{proof}
Suppose that 
\begin{equation}
  \frac{x^2-r^2-t^{-2}}{2rt^{-1}} \ge 1.
\end{equation}
Then $x^2\ge (r+t^{-1})^2$, and 
thus from \eqref{e-con3}, 
\begin{equation}
  y^2= t^2x^2 + (r^2-1)(1-t^2) \ge t^2(r+t^{-1})^2 + (r^2-1)(1-t^2)
= (r+t)^2.
\end{equation}
Hence we have 
\begin{equation}\label{e-endpf1}
  x\ge r+t^{-1}, \qquad y\ge r+t.
\end{equation}
Inserting \eqref{e-endpf1} into \eqref{e-con1}, \eqref{e-con2}, we have 
\begin{equation}
  \frac{a\gamma (1-t^2)}{r+t^{-1}} 
\ge t(1+\gamma a)\bigl( 1-\frac{r_1}{r}\bigr), 
\qquad  \frac{a(1-t^2)}{r+t}
\ge t(1+\gamma a)\bigl( \frac{r_2}{r}-1 \bigr)
\end{equation}
We multiply the first inequality by $r(r+t^{-1})$, 
and multiply the second inequality by $r(r+t)$.
Then by adding the two inequalities, we obtain, 
after some algebra, $0\ge 2$, which is a contradiction.

Now suppose that  
\begin{equation}
  \frac{x^2-r^2-t^{-2}}{2rt^{-1}} \le -1.
\end{equation}
This implies that $(tx)^2\le (1-rt)^2$.
Since $r<r_2\le \frac1{t}$ for $\gamma\ge 1$, we have 
\begin{equation}\label{e-endpf2}
  tx\le 1-rt,
\end{equation}
and from \eqref{e-con3}, 
\begin{equation}
  y^2= t^2x^2+(r^2-1)(1-t^2) 
\le (1-rt)^2+ (r^2-1)(1-t^2)  = (r-t)^2.
\end{equation}
Hence we have $y\le |r-t|$. We distinguish two cases $r\ge t$ and $r<t$.
For the first case when $r\ge t$, from \eqref{e-endpf2} we have 
\begin{equation}
  x+y\le (t^{-1}-r) + (r-t) = t^{-1}-t,
\end{equation}
which contradicts Lemma \ref{lem-endptxyt}.
For the second case when $r<t$, 
\begin{equation}\label{e-endpf3}
  y\le t-r.
\end{equation}
Inserting \eqref{e-endpf2}, \eqref{e-endpf3} 
into \eqref{e-con1}, \eqref{e-con2}, we obtain 
\begin{equation}
  \frac{a\gamma (1-t^2)}{1-rt} \le 1+\gamma a - \frac{t(a+1)}{r}, 
\qquad \frac{a(1-t^2)}{t-r} \le \frac{a+1}{r} - t(1+\gamma a).
\end{equation}
Multiply the first inequality by $(1-rt)r$, multiply 
the second inequality by $(t-r)r$, 
and add the resulting two inequalities.
Then after some algebra, we find $a\le 0$, which is a contradiction.
This proves the lemma.
\end{proof}

It now follows from the preceding discussion that 
there is a unique solution $\xi$, $Im(\xi)>0$, to 
\eqref{e-endcon}, or equivalently \eqref{e-end1}.

\subsection{The contour $\Gamma_1$ }

As in \eqref{e-z0inprop}, set 
\begin{equation}
  z_0 := \frac{a+1}{-R(0)(1+\gamma a)},
\end{equation}
where $R(0)$ is given in \eqref{e-rxyxi}.
We emphasize that $z_0$ is uniquely determined by the endpoint $\xi$, 
and is independent of the curve $\Gamma$ 
in subsection \ref{subsec-end}, as long as $\Gamma$ does not
intersect $(-\infty, 0]$,
Note that $t<z_0<t^{-1}$ from the condition $r_1<r:=-R(0) <r_2$ 
of \eqref{e-rxyineq}.

Now define 
\begin{equation}\label{e-Q}
  Q(z):= -(1+\gamma a)^2  
\frac{(z-\xi)(z-\overline{\xi})(z+z_0)^2}{z^2(z+t)^2(z+t^{-1})^2},
\end{equation}
and we consider the quadratic differential $Q(z)(dz)^2$.  In this
subsection, we are interested in the trajectories $Q(z)(dz)^2>0$.  In
the next subsection we consider the orthogonal trajectories
$Q(z)(dz)^2<0$. A general reference for quadratic differentials is
Chapter 8 of \cite{Pomm}.  Note, in particular, that if two
trajectories (or two orthogonal trajectories) meet at a point $z$
where $Q(z)$ is analytic and nonzero, then they must be identical.
The quadratic differential $Q(dz)^2$ has double poles at $0, -t,
-t^{-1}$ and $\infty$, a double zero at $-z_0$ and simple zeros at
$\xi, \overline{\xi}$. We first consider the local structure of the
trajectories near the poles and the zeros.
\subsubsection{Local structure  }

\begin{itemize}
\item Near $0$: $Q(dz)^2 \sim \frac{-c^2}{z^2}(dz)^2$, 
where $c= (1+\gamma a)|\xi|z_0 >0$.
Thus we want trajectories $z=z(t)$ satisfying 
$\frac{ic}{z}\frac{dz}{dt} \sim c_1$ for some real constant $c_1$.
The solution is given by 
$z\sim c_2e^{-i(c_1/c) t}$, and hence the trajectories near $z=0$ 
are circular.
\item Near $-t$: $Q(dz)^2 \sim \frac{-c^2}{(z+t)^2}(dz)^2$, 
where $c= \frac{(1+\gamma a)|t+\xi||-t+z_0|}{t(t^{-1}-t)} >0$. 
As in the $z=0$ case, the trajectories are circular.
\item Near $-t^{-1}$: $Q(dz)^2 \sim \frac{-c^2}{(z+t^{-1})^2}(dz)^2$,
where $c= \frac{(1+\gamma a)|t^{-1}+\xi||-t^{-1}+z_0|}{t^{-1}(t^{-1}-t)} >0$.
Again, the trajectories are circular.
\item Near $\infty$, $Q(dz)^2 \sim 
\frac{-c^2}{z^2}(dz)^2= \frac{-c^2}{w^2}(dw)^2$,
where $c= 1+\gamma a >0$, $w=\frac1z$.
Thus once again the trajectories are circular.
\item Near $-z_0$: $Q(dz)^2 \sim -c^2 (z+z_0)^2(dz)^2$, 
where $c=\frac{(1+\gamma a)|z+\xi|}{z_0(z_0-t)(t^{-1}-z_0)}>0$.
Thus we seek trajectories $z=z(t)$ satisfying 
$ic(z+z_0)\frac{dz}{dt} \sim c_1$ for some real constant $c_1$.
The solution with $z(0)=-z_0$ is given by $(z+z_0)^2 \sim -i(2c_1/c)t$, thus 
we have $arg(z+z_0)= \frac\pi4+ \frac{k\pi}4$, $k=0,1,2,3$.
Hence there are 4 trajectories starting from $-z_0$, all 
making an angle $\frac\pi4$ with the real line.
\item Near $\xi$: $Q(dz)^2 \sim c(z-\xi)(dz)^2$ for some $c\in\C$. 
Thus we seek trajectories $z=z(t)$ satisfying
$c\sqrt{z-\xi}\frac{dz}{dt} \sim c_1$ for some real constant $c_1$.
The solution with $z(0)=\xi$ is given by $(z-\xi)^{3/2} \sim (c_1/c)t$, 
and hence $arg(z-\xi)= c_2+ \frac{2k\pi}3$, $k=0,1,2$.
Thus there are 3 trajectories starting from $\xi$, 
making an angle $\frac{2\pi}3$ between themselves.
\item Near $\overline{\xi}$: $Q(dz)^2 \sim c(z-\overline{\xi})(dz)^2$ 
for some $c\in\C$.
Again,
there are 3 trajectories starting from $\overline{\xi}$, 
making an angle $\frac{2\pi}3$ between themselves.
\end{itemize}

\begin{figure}[ht]
 \centerline{\epsfig{file=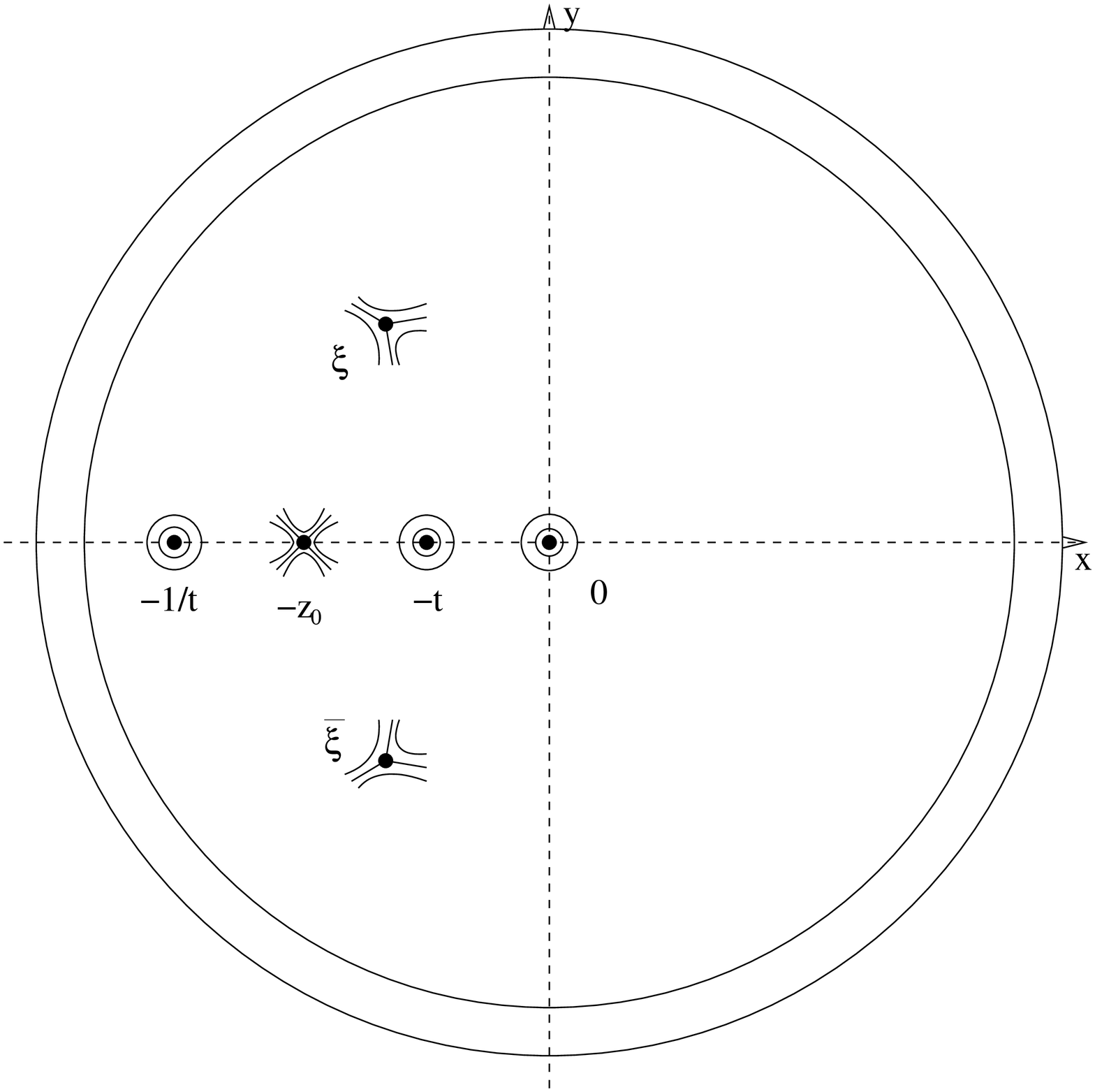, width=9cm}}
 \caption{Local structure of the trajectories of $Q(dz)^2$ near the poles 
and the zeros.}
\label{fig-local}
\end{figure}

The local structure of the trajectories of $Q(dz)^2$ 
near the poles and the zeros are summarized in 
Figure \ref{fig-local}.

\subsubsection{Global structure }

Note that as $Q(z)=\overline{Q(\overline{z})}$, 
if $\{z(t): \alpha < t<\beta\}$ is a trajectory of $Q(dz)^2>0$, 
then $\{\overline{z(t)}: \alpha < t<\beta\}$ is also a trajectory.
Also from \eqref{e-Qdef1}, we see that $Q(z)<0$ for 
$z\in\R\setminus\{-z_0, 0, -t, -t^{-1}\}$,
hence all trajectories that cross the real axis do so at $\pi/2$,
and if $\{z(t): \alpha<t<\beta\}$ is a trajectory that 
satisfies $Im(z(t))>0$ for $\alpha<t<\beta$, 
and $Im(z(\beta))=0$, $Re(z(\beta))\neq -t^{-1}, -z_0, -t, 0$, then 
$\{ \overline{z(2\beta-t)}: \beta < t<2\beta-\alpha\}$ 
gives a smooth continuation 
of $\{z(t)\}$ into $Im(z)<0$.

We need the following lemma.

\begin{figure}[ht]
 \centerline{\epsfig{file=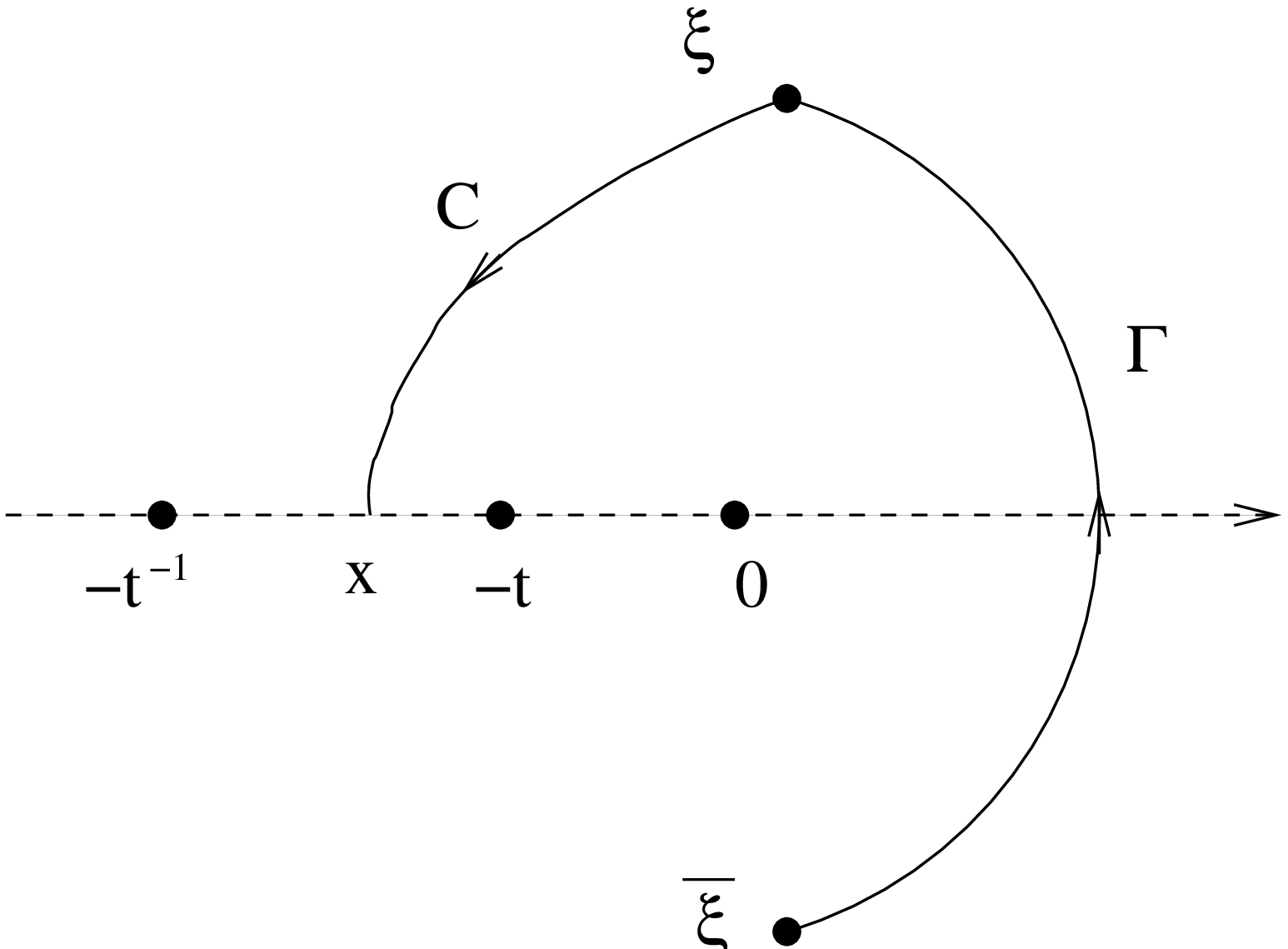, width=6cm}
\epsfig{file=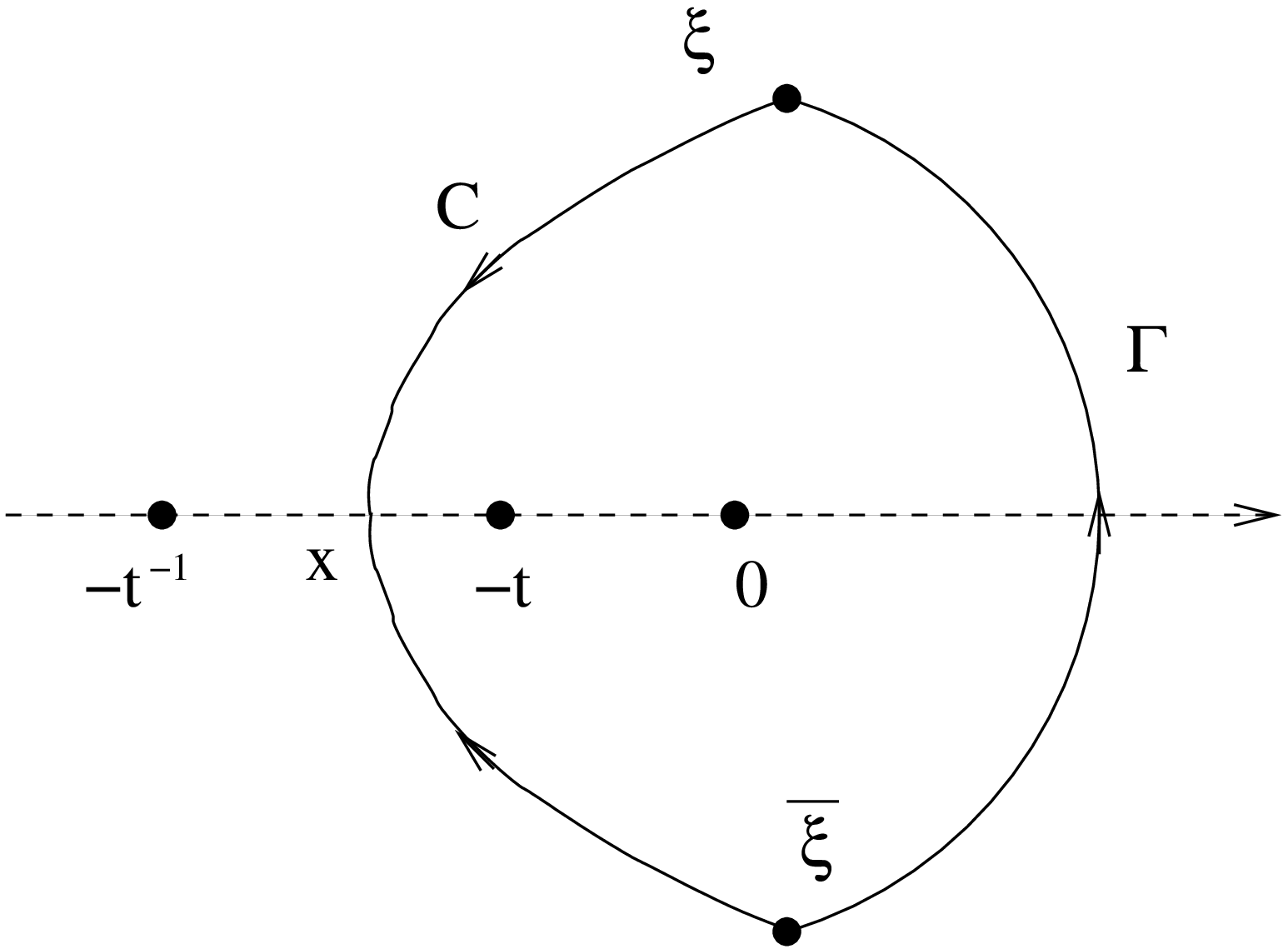, width=6cm}}
 \caption{The contours $C$ and $C'$.}
\label{fig-Ccontour}
\end{figure}

\begin{lem}\label{lem-z0}
Let $\xi, \overline{\xi}$ be as in subsection \ref{subsec-end}.
Let $\Gamma$ be a simple curve with endpoints $\xi, \overline{\xi}$ 
which does not intersect 
$(-\infty, 0]$ and is symmetric under reflection about the real axis.
Choose the branch of $\sqrt{Q(z)}$ to be 
analytic in $\C\setminus\overline{\Gamma}$, 
and satisfy $\sqrt{Q(z)} \sim \frac{i(1+\gamma a)}{z}$ as $z\to\infty$.
For any real number $x$ satisfying $-t^{-1}<x<-t$, 
let $C$ be a smooth curve in $\C_+$ with endpoints $\xi$ and $x$,
oriented from $\xi$ to $x$, which does not intersect
$\overline{\Gamma}$ (see Figure \ref{fig-Ccontour}).
Then we have 
\begin{equation}\label{e-basiclem}
  Re \int_{C} \sqrt{Q} dz =0.
\end{equation}
\end{lem}

\begin{proof}
  Let $C'$ be the closure of $C\cup\overline{C}$, 
oriented from $\xi$ to $\overline{\xi}$. 
Hence $C'$ is a curve which has endpoints $\xi, \overline{\xi}$, 
intersects the real axis at $-t^{-1}<x<-t$, and satisfies
$C'=\overline{C'}$ (see Figure \ref{fig-Ccontour}).
Let $C^*=\{\overline{z}: z\in C\}$, oriented from $\overline{\xi}$ to $x$.
Using $\sqrt{Q(z)}=-\overline{\sqrt{Q(\overline{z})}}$ 
and the realness of $x$, we have 
\begin{equation}
\begin{split}
  Re \int_{C} \sqrt{Q(z)}dz &= 
\frac1{2} \biggl[ 
\int_{C} \sqrt{Q(z)}dz + \overline{\int_{C} \sqrt{Q(z)}dz} \biggr]
= \frac1{2} \biggl[ \int_{C} \sqrt{Q(z)}dz + 
\int_{C^*} \overline{\sqrt{Q(\overline{z})}}dz \biggr]  \\
& = \frac1{2} \biggl[ \int_{C} \sqrt{Q(z)}dz +
\int_{-C^*} \sqrt{Q(z)}dz \biggr] 
= \frac1{2} \int_{C'} \sqrt{Q(z)} dz.
\end{split}
\end{equation}
Hence we want to prove that the last integral is $0$.

Set (cf. \eqref{e-h}, \eqref{e-hres-0})
\begin{equation}\label{e-hset0}
\begin{split}
  h(z) &= \frac{R(z)}{2\pi i} \int_{\Gamma} \frac{W'(s)}{R_+(s)(s-z)} ds \\
&= \frac12 W'(z) + \frac12 R(z)\biggl(
\frac{\gamma a}{(z+t^{-1})R(-t^{-1})} + \frac{a}{(z+t)R(-t)}
-\frac{a+1}{zR(0)} \biggr) \\
&= \frac12 W'(z) + \frac12(1+\gamma a)
\frac{R(z)(z+z_0)}{z(z+t)(z+t^{-1})},
\end{split}
\end{equation}
where the second equality follows from a residue calculation, 
while the third equality follows from 
the endpoint conditions \eqref{e-endcon}.
Thus from the definition \eqref{e-Q} of $Q$ 
and the choice of $\sqrt{Q}$, we have 
$\sqrt{Q} = i(2h-W')$.
Now 
\begin{equation}
\begin{split}
  \frac1{i}\int_{C'} \sqrt{Q(z)}dz &= \int_{C'} (2h-W')dz 
= \int_{C'} \frac{R(z)}{\pi i} dz \int_{\Gamma} \frac{W'(s)}{R_+(s)(s-z)}ds
- \int_{C'} W'(s)ds \\
& = \frac1{\pi i} \int_{\Gamma} \frac{W'(s)}{R_+(s)} ds 
\int_{C'} \frac{R(z)}{s-z}dz - \int_{C'} W'(s)ds.
\end{split}
\end{equation}
By a residue calculation, for $s\in \Gamma$, 
\begin{equation}
  \int_{C'} \frac{R(z)}{s-z}dz = \lim_{r\to\infty} \frac12 
\int_{|z|=r} \frac{R(z)}{s-z} dz - \pi i R_+(s) 
= \pi i \biggl( \frac{\xi+\overline{\xi}}2 -s - R_+(s) \biggr).
\end{equation}
Hence using the endpoints conditions \eqref{e-end1}, we have 
\begin{equation}
\begin{split}
  \frac1{i}\int_{C'} \sqrt{Q(z)}dz
& = \int_{\Gamma} \frac{W'(s)}{R_+(s)} 
\biggl( \frac{\xi+\overline{\xi}}2 -s - R_+(s) \biggr) ds 
- \int_{C'} W'(s)ds \\
& = 2\pi i - \int_{\Gamma} W'(s)ds - \int_{C'} W'(s)ds \\
& = 2\pi i - \int_{\Gamma'} W'(s)ds,
\end{split}
\end{equation}
where $\Gamma'$, the closure of $\overline{ \Gamma\cup C'}$, 
encloses $-t, 0$, but not $-t^{-1}$, 
and is oriented counter-clockwise.
But a direct calculation shows 
$\int_{\Gamma'} W'(s)ds = 2\pi i$, and we obtain the lemma.
\end{proof}

\begin{rem}
The authors are indebted to Nick Ercolani who suggested
that a formula such as \eqref{e-basiclem} should be true.
\end{rem}

\begin{figure}[ht]
 \centerline{\epsfig{file=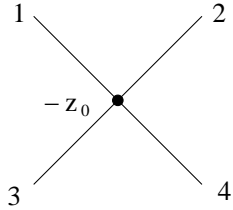, width=3cm}}
 \caption{Local structure of the trajectories of $Q(dz)^2$ near $-z_0$.}
\label{fig-z0local}
\end{figure}

The following general results for the trajectories 
of quadratic differentials 
are given in Lemmas 8.3, 8.4 of \cite{Pomm}, and  
are used in several places in the arguments that follow.

\begin{lem}\label{lem-Pomm}
Let $Q(z)(dz)^2$ be a quadratic differential 
in a simply connected domain $G$.
\begin{enumerate}
\item 
If there is at most one pole of $Q(z)$ in $G$ and this pole is simple,  
then there is no closed Jordan curve in $G$ 
consisting only of trajectories (or orthogonal trajectories) 
and their endpoints.
\item 
Suppose that $Q(z)$ has no poles in $G$
and let $\Gamma$ be a trajectory (or an orthogonal trajectory). 
Then in both directions, $\Gamma$ ends at a zero of $Q$ 
or converges to $\partial G$.
\end{enumerate}
\end{lem}

\begin{figure}[ht]
 \centerline{\epsfig{file=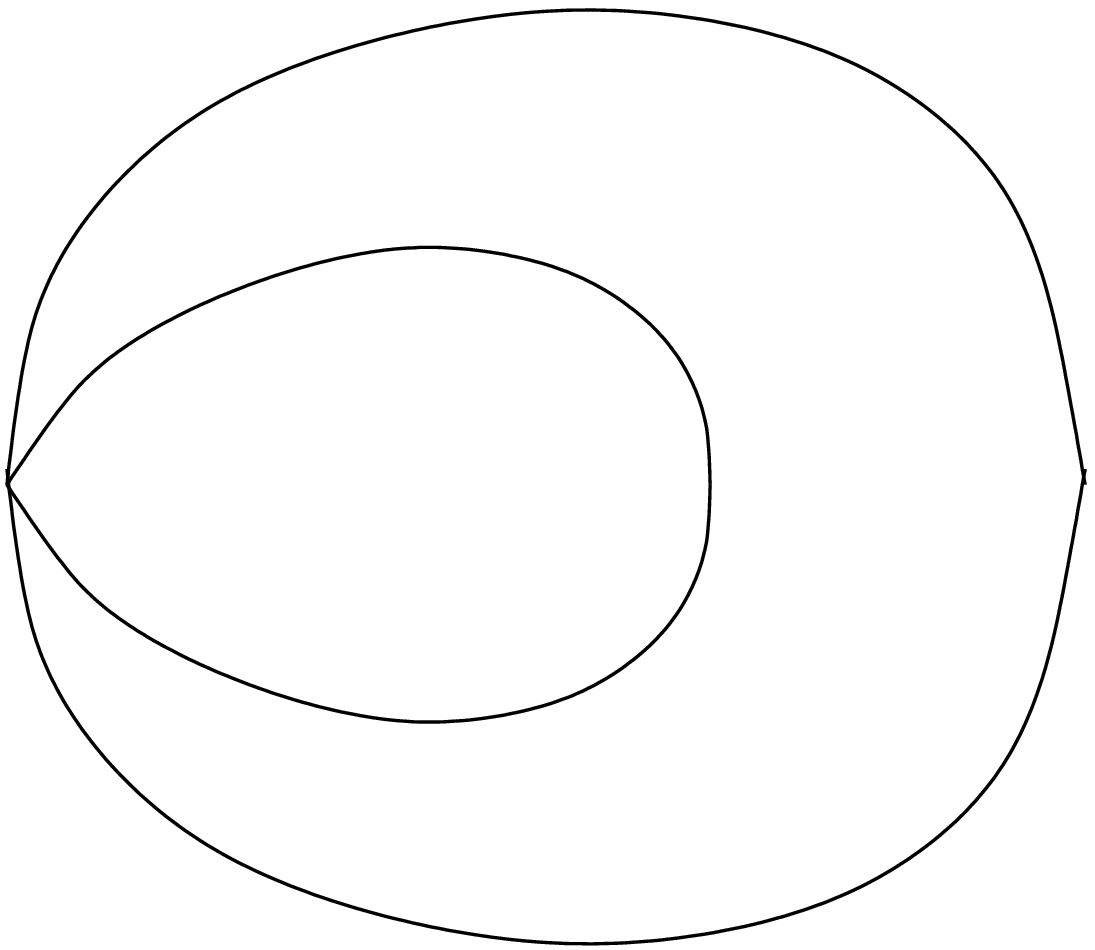, width=3cm} \qquad\qquad
\epsfig{file=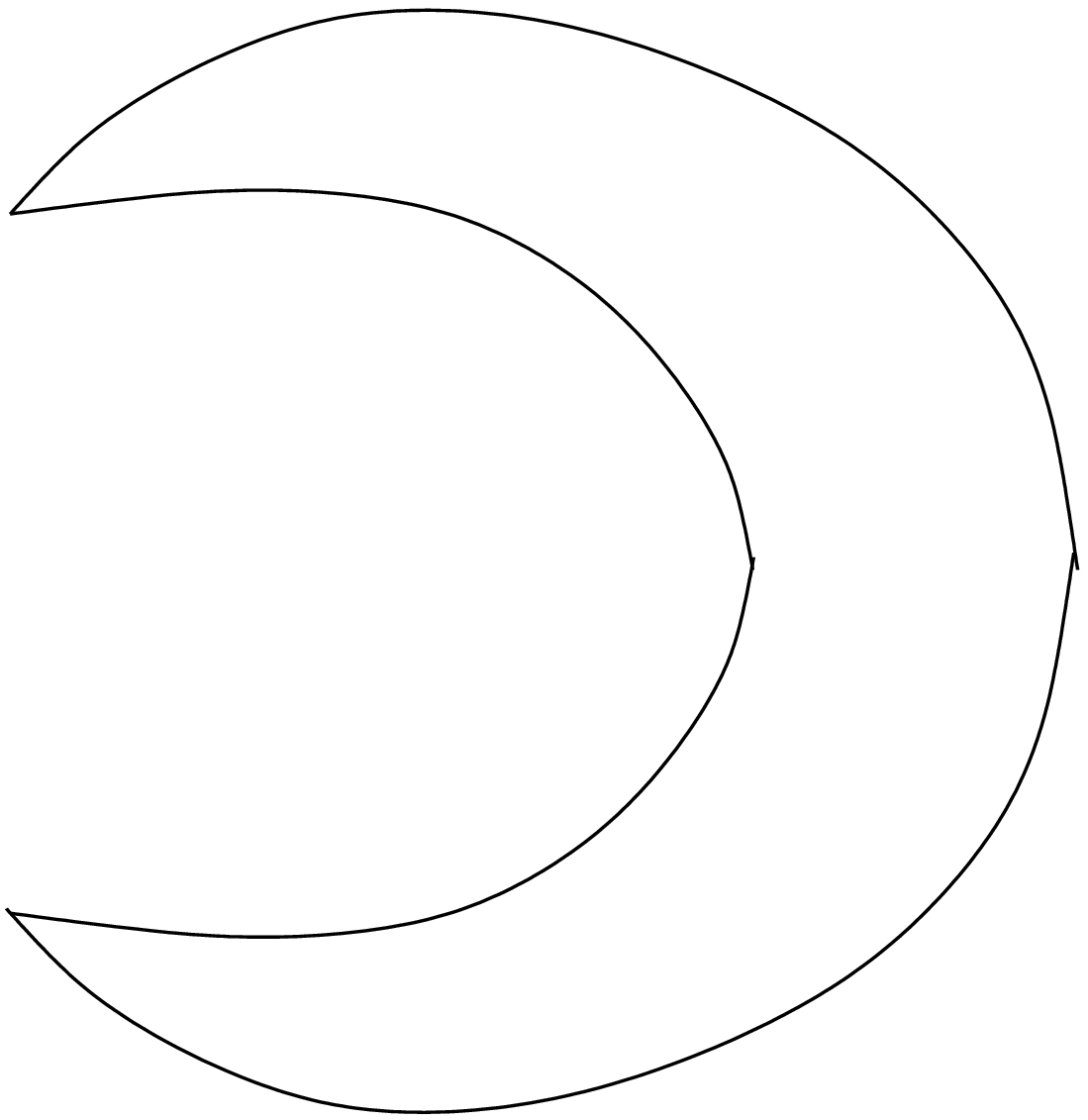, width=3cm}}
 \caption{Unfolding of a ``closed'' trajectory meeting itself at a point.}
\label{fig-lem}
\end{figure}

\begin{rem}
At various points in the argument that follows, 
we will use Lemma \ref{lem-Pomm} (i) in a slightly extended form.
For example, we will want to consider 
Jordan curves consisting of trajectories of the form in 
the first picture of Figure \ref{fig-lem}
for which the hypotheses of the Lemma are not fully satisfied.
However, if we, for example, make a change of 
variables $z\mapsto \zeta=z^{1-\epsilon}$, 
$0<\epsilon<1$, then the figure takes the form 
as the second picture in Figure \ref{fig-lem}
and as the change of variable takes trajectories to trajectories, 
it is easy to verify that the hypotheses of Lemma \ref{lem-Pomm} (i) 
are now satisfied.
We will use this extended form of the Lemma without further comment below.
\end{rem}

Denote the four trajectories emerging from $-z_0$ by $1,2,3,4$ 
as shown in Figure \ref{fig-z0local}.
First, consider the trajectory emerging from $-z_0$ along the ray $1$.
Since there are no poles in $\C_+$, which is a simply connected region, 
this trajectory must either go to 
$\xi$, the zero of $Q$ in $\C_+$, or escape from $\C_+$.
Suppose that the ray $1$ does not escape from $\C_+$.
Then from the definition of a trajectory, 
$\int_1 \sqrt{Q}dz \in \R\setminus\{0\}$, 
and hence the ray $1$ (and similarly the ray $2$) 
can not go to $\xi$ due to Lemma \ref{lem-z0}. 
So the ray $1$ must exit from $\C_+$. 
Now it can not exit through $-z_0$ 
because then by the local structure of the trajectories,
the ray $1$ comes back to $z_0$ through ray $2$, and the 
interior of the loop has no poles of $Q$, 
contradicting Lemma \ref{lem-Pomm} (i).
Also, again by the local structure of the trajectories, the ray 1 can not exit 
through $-t^{-1}, -t, 0$ or $\infty$. 
Now there are five possibilities: 
the ray 1 exits from $C_+$ through $(-\infty, -t^{-1})$, 
$(-t^{-1}, -z_0)$, $(-z_0, -t)$, $(-t, 0)$, or $(0,\infty)$.
We examine each case.

\begin{itemize}
\item[(i)]
The ray 1 can not exit through $(-t^{-1}, -z_0)$,
for if it does, the trajectory 
can be continued by complex conjugation as remarked before, 
and the interior of the loop contains no poles of $Q(dz)^2$,
contradicting Lemma \ref{lem-Pomm} (i).
\item[(ii)]
By the same argument, the ray can not exit through $(-z_0, -t)$.
\item[(iii)] 
Suppose the ray 1 exits at $z_1\in (-t, 0)$. 
Then extending the ray by conjugation, 
we obtain the closed loop 
(see the first Picture in Figure \ref{fig-Gamma1_1,2}). 
Now the trajectory along ray 2 from $-z_0$ also can not go to $\xi$ 
and hence must exit $\C_+$. 
\begin{figure}[ht]
 \centerline{\epsfig{file=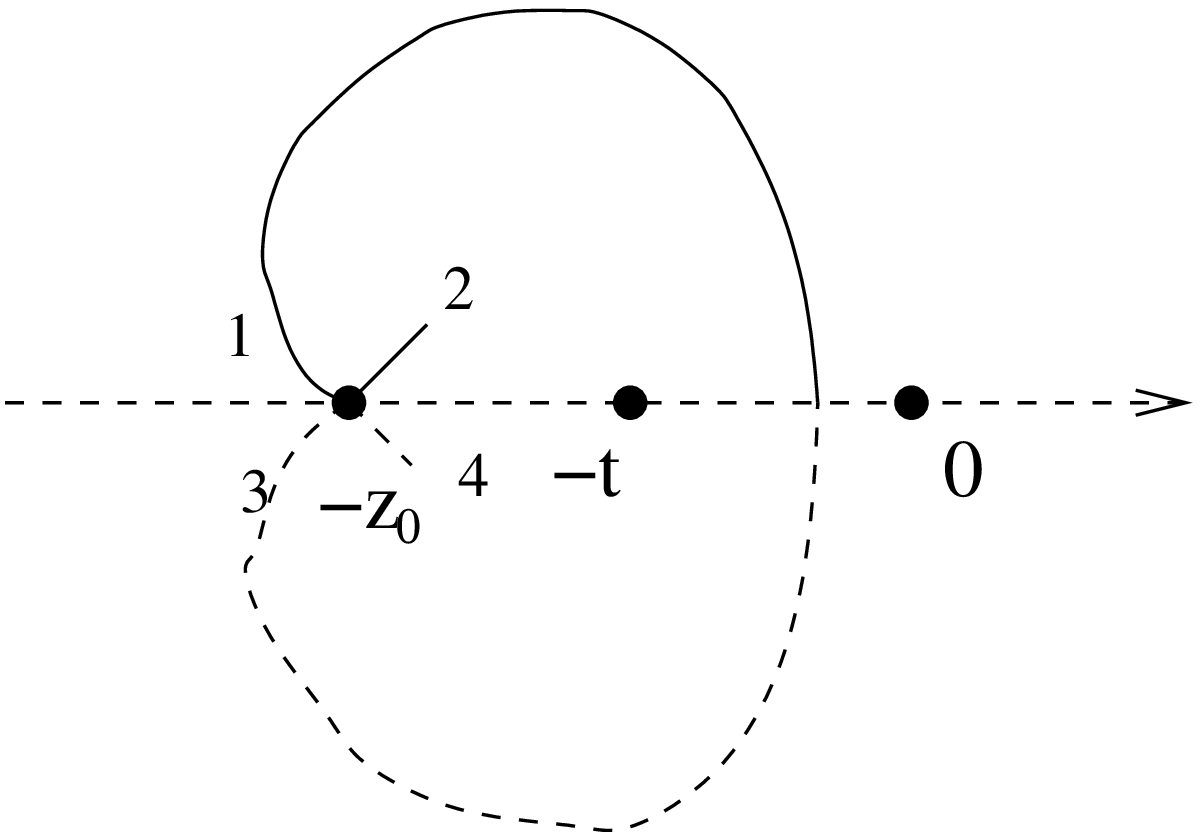, width=6cm}
\epsfig{file=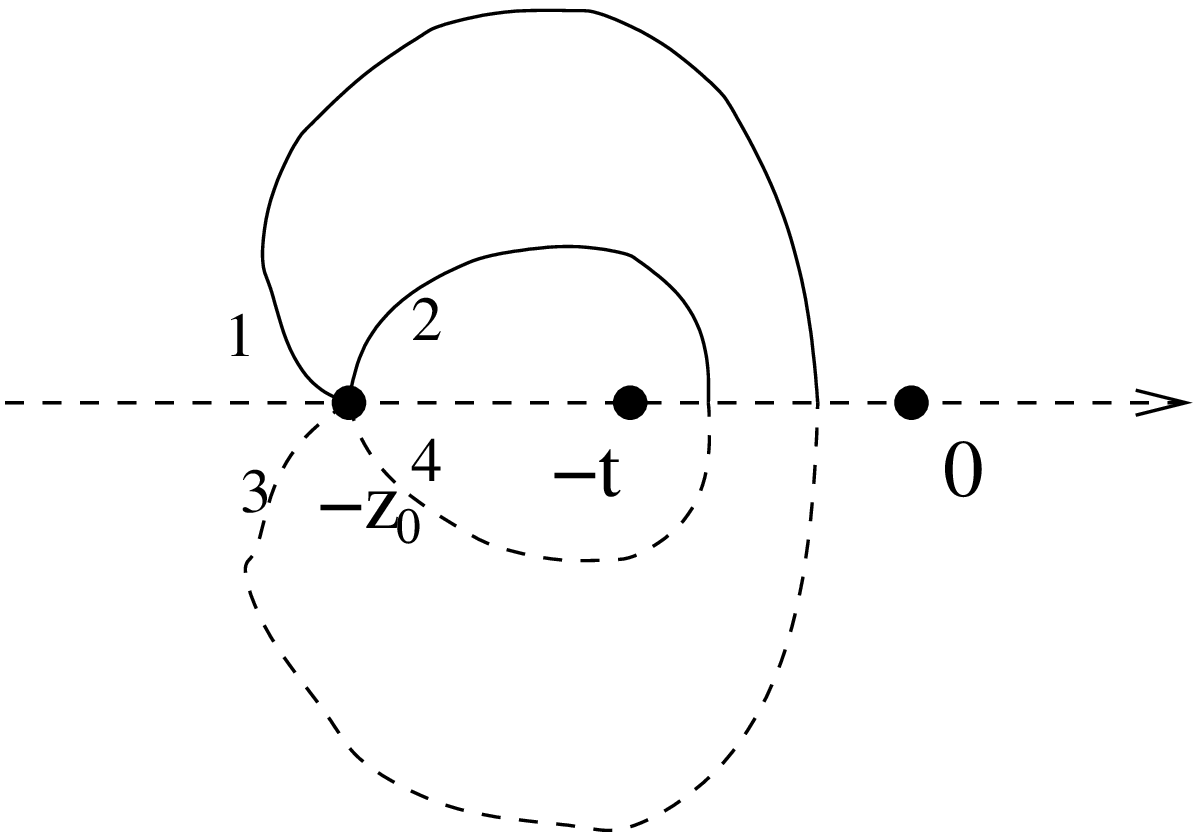, width=6cm}}
 \caption{Ray 1, case (iii).}
\label{fig-Gamma1_1,2}
\end{figure}
Also the ray 2 can not cross the trajectory emerging from the ray 1.
Hence it must cross the real axis between $(-z_0,-t)$ or $(-t, z_1)$.
The first case can not hold as the closed loop obtained by 
extending the ray 2 by conjugation has no poles inside.
In the latter case, the simply connected region 
formed by the two loops has no poles, which is again a contradiction 
(see the second picture in Figure \ref{fig-Gamma1_1,2}).
Thus ray 1 can not exit through $(-t, 0)$.
\item[(iv)] Suppose that the ray 1 exits at $z_2\in (0,\infty)$.
Then by arguing as in case (iii), the trajectory along the ray 2 
from $-z_0$ must exit $\C_+$ as some point $z_3\in (-t, 0)$. 
\begin{figure}[ht]
 \centerline{\epsfig{file=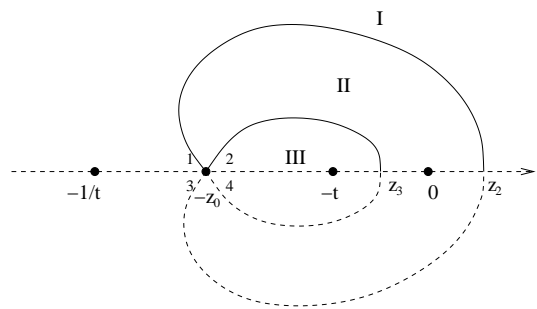, width=7cm}}
 \caption{Ray 1, case (iv).}
\label{fig-Gamma1_3}
\end{figure}
We denote the open regions in $\C_+$ divided by the trajectories emerging from 
$-z_0$ by I, II, III as shown in Figure \ref{fig-Gamma1_3}.
Now the zero $\xi$ of $Q$ lies in one of the regions I, II, III.
(Note that $\xi$ can not be on the trajectory emerging from $-z_0$ 
due to Lemma \ref{lem-z0}.)
\begin{itemize}
\item[(iv-1)]
If $\xi\in III$, then as before, 
each of the three trajectories emerging from $\xi$ exit 
through $(-z_0, -t)$ or $(-t, z_3)$.
Hence at least two of the trajectories exit through the same interval.
Then the simply connected region bounded by these two trajectories 
has no pole, which is again a contradiction.
\item[(iv-2)]
If $\xi\in II$, we reach a contradiction by a similar argument.
\item[(iv-3)]
If $\xi\in I$, the three trajectories emerging from $\xi$ 
must exit, one through $(-\infty, -t^{-1})$, one through 
$z_4\in (-t^{-1}, -z_0)$, 
and one through $(z_2, \infty)$.
But this in turn gives a contradiction, 
because $Re \int_{\xi}^{z_4} \sqrt{Q}dz=0$ by Lemma \ref{lem-z0}.
\end{itemize}
Hence in each case, we have a contradiction.
\end{itemize}
Thus we conclude that the trajectory along the ray 1 emerging from $-z_0$
exits $\C_+$ through $z_5\in (-\infty, -t^{-1})$.

Now we consider the trajectory along the ray 2 emerging from $-z_0$.
As before, it must exit $\C_+$ through either of 
$(-\infty, z_5)$, $(-z_0, -t)$, $(-t, 0)$ or $(0,\infty)$.
\begin{itemize}
\item[(i)]
If it exits through $(-\infty, z_5)$,  
the two trajectories from the rays 1 and 2, 
extended as in the remark above to $\C_+$, form a simply connected region 
which does not contain poles, this is again a contradiction.
\item[(ii)]
Suppose that the  
trajectory along the ray 2 exits through $(-z_0, -t)$.
Then the ($\C_-$-extended) loop of the trajectory 2 contains no poles, 
and again we have a contradiction.
\item[(iii)]
\begin{figure}[ht]
 \centerline{\epsfig{file=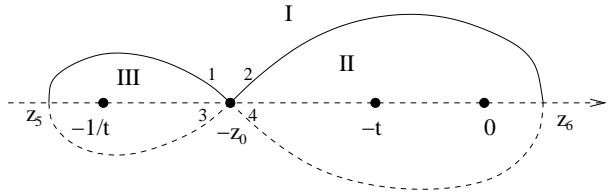, width=8cm}}
 \caption{Ray 2, case (iii).}
\label{fig-Gamma1_4}
\end{figure}
Suppose that the trajectory along the ray 2 exits through
$z_6\in (0,\infty)$.
Then $\xi$ lies either in the region I, II or III 
as of Figure \ref{fig-Gamma1_4}.
If $\xi \in I$, at least two of the trajectories from $\xi$ exit together 
through $(-\infty, z_5)$ or $(z_6,\infty)$.
This yields a contradiction as in (iv-1) above.
The case $\xi\in III$ leads to a similar contradiction.
If $\xi\in II$, then the three trajectories exit, one through 
$(-z_0, -t)$, one through $(-t, 0)$, and one through $(0, z_6)$.
But this gives in turn a contradiction by Lemma \ref{lem-z0} 
as in the case (iv-3) above.
\end{itemize}
Therefore the trajectory along the ray 2 emerging from $-z_0$ 
must exit $\C_+$ through $z_7\in (-t, 0)$.

\begin{figure}[ht]
 \centerline{\epsfig{file=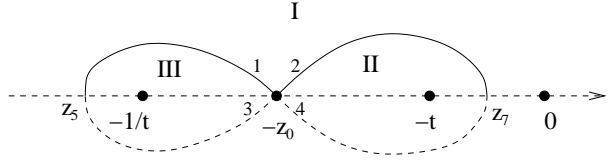, width=8cm}}
 \caption{Rays 1 and 2.}
\label{fig-Gamma1_5}
\end{figure}

\begin{figure}[ht]
 \centerline{\epsfig{file=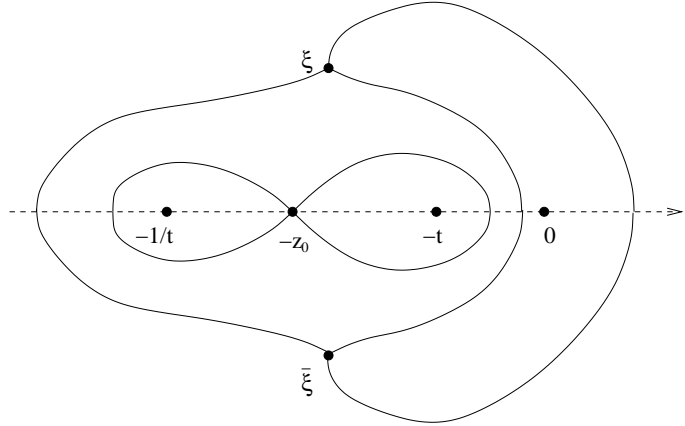, width=9cm}}
 \caption{Global structure of the trajectories of $Q(dz)^2$.}
\label{fig-Gammaglobal}
\end{figure}

Thus $\xi$ lies either in regions I, II or III of Figure \ref{fig-Gamma1_5}.
But by a now familiar argument as above, 
$\xi$ can not lie in II or III. 
Hence $\xi\in I$.
Then the three trajectories emerging from $\xi$ exit at
some points $z_8, z_9, z_{10}$ 
where $z_8\in (-\infty, z_5), z_9\in (z_7, 0), z_{10}\in (0,\infty)$. 
This shows that the global structure of the trajectories of $Q(dz)^2$ 
is given in Figure \ref{fig-Gammaglobal}.  

The above considerations prove, in particular, 
that there is a trajectory emerging from 
$\xi$ and ending at $\overline{\xi}$ 
which crosses the real axis to the right of $0$.
We take $\Gamma_1$ to be this curve, oriented 
from $\overline{\xi}$ to $\xi$.
Define $h$ by \eqref{e-hset0} 
with the choice of the contour $\Gamma=\Gamma_1$: 
let $R(z)=\sqrt{(z-\xi)(z-\overline{\xi})}$ be analytic 
in $\C\setminus\Gamma_1$, $R(z)\sim z$ as $z\to\infty$. 
Thus 
\begin{equation}\label{e-hdef}
\begin{split}
  h(z) &= \frac{R(z)}{2\pi i} \int_{\Gamma_1} \frac{W'(s)}{R_+(s)(s-z)} ds \\
&= \frac12 W'(z) + \frac12 R(z)\biggl(
\frac{\gamma a}{(z+t^{-1})R(-t^{-1})} + \frac{a}{(z+t)R(-t)}
-\frac{a+1}{zR(0)} \biggr) \\
&= \frac12 W'(z) + \frac12 (1+\gamma a)
\frac{R(z)(z+z_0)}{z(z+t)(z+t^{-1})}.
\end{split}
\end{equation}
Condition (a) of Proposition \ref{prop-h}, 
$h_+(z)+h_-(z)=W'(z)$ for $z\in\Gamma_1$, 
now follows by properties of the Cauchy operator, 
and the condition (b) by the endpoint condition \eqref{e-end1}.

Now we consider condition (c).
Define $\sqrt{Q(z)}=i(1+\gamma a)
\frac{R(z)(z+z_0)}{z(z+t)(z+t^{-1})}$, so that it is analytic 
in $\C\setminus \overline{\Gamma_1}$, 
and $\sqrt{Q} \sim \frac{i(1+\gamma a)}{z}$ as $z\to\infty$.
Then 
\begin{equation}\label{e-hdiffQ}
  i(h_+(z)-h_-(z))=\sqrt{Q(z)}_+, \qquad z\in\Gamma_1.
\end{equation}
Along the trajectory $\Gamma_1$, we must have 
$\sqrt{Q}dz\in \R\setminus\{0\}$.
Let $z_{10}\in(0,\infty)$ be the point at which $\Gamma_1$ crosses 
the real axis.
As $R(z)=\sqrt{(z-\xi)(z-\overline{\xi})}>0$ 
for $z\in (z_{10},\infty)$, it follows that 
for all upward-oriented trajectories that cross the real axis at points 
$x\in (z_{10},\infty)$,  $\sqrt{Q}dz|_{z=x}<0$.
It particular, for $\Gamma_1$, $\sqrt{Q}_-dz|_{z_{10}}<0$, 
and hence $\sqrt{Q(z)}_-dz<0$ for all $z\in\Gamma_1$.
Now since $\sqrt{Q}_+=-\sqrt{Q}_-$ on $\Gamma_1$, 
we have 
\begin{equation}
  \sqrt{Q}_+dz \in \R_+, \quad \text{on $\Gamma_1$,}
\end{equation}
and hence from \eqref{e-hdiffQ},
\begin{equation}
  i(h_+(z)-h_-(z))dz >0, \quad z\in\Gamma_1,
\end{equation}
which proves condition (c).

\subsection{The contour $\Gamma_2$}

\begin{figure}[ht]
 \centerline{\epsfig{file=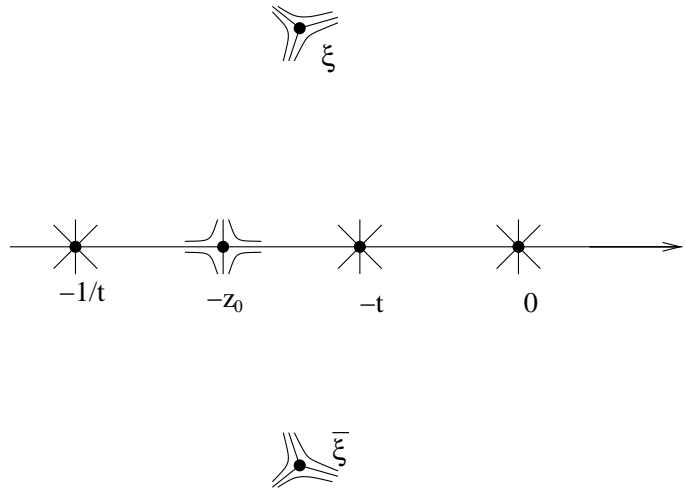, width=9cm}}
 \caption{Local structure of the orthogonal trajectories 
near the finite poles and zeros of $Q(dz)^2$.}
\label{fig-Gamma2local}
\end{figure}

Again we choose the branch of $\sqrt{Q}$ so that $\sqrt{Q}$ is analytic 
in $\C\setminus\Gamma_1$ and $\sqrt{Q} \sim \frac{i(1+\gamma a)}{z}$
as $z\to\infty$.
Now we consider the orthogonal trajectories $Q(z)(dz)^2<0$.
As before, the local structure is easy to determine. 
We summarize the local structure of the orthogonal trajectories 
near the finite poles and the zeros of $Q(dz)^2$ in 
Figure \ref{fig-Gamma2local}.
Near $\infty$, any straight rays emerging from 
$\infty$ are orthogonal trajectories.
We note that the real axis is an orthogonal trajectory.
Hence the orthogonal trajectories can cross the real axis only 
at $0, -t, z_0, -t^{-1}$ or $\infty$.

\begin{figure}[ht]
 \centerline{\epsfig{file=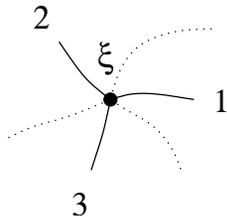, width=3cm}}
 \caption{Orthogonal trajectories and trajectories emerging from $\xi$.}
\label{fig-Gamma2_0}
\end{figure}

\begin{figure}[ht]
 \centerline{\epsfig{file=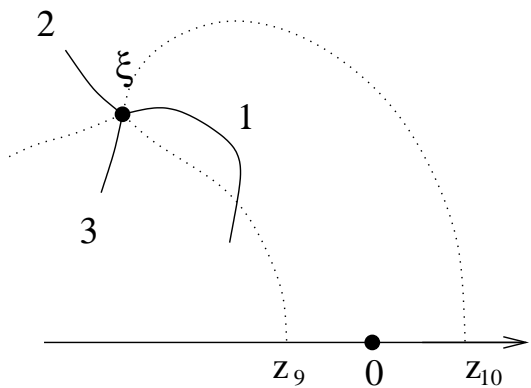, width=7cm}}
 \caption{Orthogonal trajectories emerging from $\xi$.}
\label{fig-Gamma2_1}
\end{figure}

\begin{figure}[ht]
 \centerline{\epsfig{file=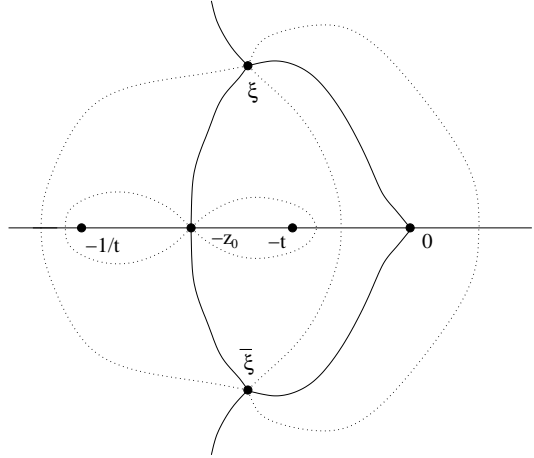, width=7cm}}
 \caption{Global structure for the orthogonal trajectories 
and the trajectories emerging from 
the zeros $\xi$, $\overline{\xi}$ and $-z_0$.}
\label{fig-Gamma2global}
\end{figure}

Now we consider the global structure.
Again by the symmetry $Q(z)=\overline{Q(\overline{z})}$, 
the orthogonal trajectories are symmetric under reflection 
about the real axis.
There are three orthogonal trajectories, denoted by $1,2,3$,  
emerging from $\xi$, each of which bisects the angle between 
two trajectories emerging from $\xi$
(see Figure \ref{fig-Gamma2_0}).
Now we show that 
the orthogonal trajectory $1$ can not cross either of the two adjacent 
trajectories emerging from $\xi$.
If so, there is a closed loop of the form 
shown the Figure \ref{fig-Gamma2_1} 
consisting of a part of a trajectory and a part of an orthogonal trajectory, 
and by analyticity of $\sqrt{Q}$, the integral 
of $\sqrt{Q}dz$ is zero around the loop.
But the integral of $\sqrt{Q}dz$ along the trajectory
is in $\R\setminus\{0\}$, and along the orthogonal trajectory, the integral 
of $\sqrt{Q}dz$ along the orthogonal trajectory is in $i\R\setminus\{0\}$.
Hence the sum can not be zero, which is a contradiction.
Therefore the orthogonal 
trajectory $1$ must exit $\C_+$ between $z_9$ and $z_{10}$.
But then from the local structure, it must exit at $0$.
Similarly, the orthogonal trajectory $2$ must go to $\infty$.
The orthogonal trajectory $3$, by a similar argument, must exit $\C_+$
either through $-t^{-1}, z_0$ or $-t$.
Suppose it exits through $-t$.
By the local structure, the orthogonal trajectory approaches along 
an angle; indeed it is easy to show that 
$z(s) \sim -t + e^{i\phi}e^{-cs}$ as $s\to\infty$
for some $0<\phi<\pi$, $c>0$.
Also by the definition of an orthogonal trajectory, we have 
\begin{equation}\label{e-Gamma2.1}
  \lim_{s\to\infty} Re \int_{\xi}^{z(s)} \sqrt{Q(z)}dz =0
\end{equation}
along the orthogonal trajectory $3$.
For $s$ large, but fixed, write 
$z(s):=-t+\epsilon e^{i\theta}$ for $\epsilon>0$, $0<\theta<\pi$,
and consider the curve $C_0$ from $\xi$ to $z(s)$
along the orthogonal trajectory $3$.
Let $C_s$ be the curve 
$\{ -t + \epsilon e^{i\beta}: \theta\le \beta\le \pi\}$, 
oriented from $z(s)$ to $-t-\epsilon$.
Thus $C_0\cup C_s$ is a curve from $\xi$ to a point in $(-t^{-1}, -t)$.
Then by Lemma \ref{lem-z0}, we have 
$\int_{C_0\cup C_s} Re \sqrt{Q}dz=0$. 
Hence we have 
\begin{equation}
\begin{split}
  \lim_{s\to\infty} Re \int_{C_0} \sqrt{Q}dz 
&= \lim_{s\to\infty} Re \int_{C_0 \cup C_s} \sqrt{Q} dz 
- Re \int_{C_s} \sqrt{Q}dz\\
&= - \lim_{s\to\infty} Re \int_{C_s} \sqrt{Q}dz  \\
&= - \lim_{s\to\infty} Re \int_{C_s} i(1+\gamma a)
\frac{R(z)(z+z_0)}{z(z+t^{-1})} 
\frac{dz}{z+t} \\
&= (1+\gamma a)\frac{R(-t)(-t+z_0)}{-t(-t+t^{-1})} (\pi -\theta) >0.
\end{split}
\end{equation}
This contradicts \eqref{e-Gamma2.1}, and hence the trajectory $3$
cannot exit $\C_+$ through $-t$.
Similar argument shows that it cannot exit $\C_+$ through $-t^{-1}$.
Hence the trajectory $3$ exits $\C_+$ through $-z_0$.
Thus we have obtained the global structure for 
the orthogonal trajectories of 
$Q(dz)^2$ emerging from the zeros $\xi$, $\overline{\xi}$, $-z_0$
as shown by the solid curves in the Figure \ref{fig-Gamma2global}.
The dotted curves in the Figure \ref{fig-Gamma2global} 
denote the trajectories 
already displayed in Figure \ref{fig-Gammaglobal}.

In particular, we have shown 
that there is an orthogonal trajectory
(more precisely, a union of two orthogonal trajectories, 
one in $\C_+$ and the other, its conjugate, in $\C_-$, 
meeting at the point $-z_0$ which is a zero of $Q(dz)^2$)
emerging from $\xi$
and ending at $\overline{\xi}$ and which crosses the real axis at $-z_0$. 
We denote the curve by $\Gamma_2$, and take the orientation 
from $\xi$ to $\overline{\xi}$.
Since $\Gamma_2$ is a union of two orthogonal trajectories, 
we have $\sqrt{Q}dz\in i\R\setminus\{0\}$ for $z\in\Gamma_2\setminus\{0\}$
and an explicit computation using 
$\sqrt{Q(z)}=i(1+\gamma a) \frac{R(z)(z+z_0)}{z(z+t)(z+t^{-1})}$ shows 
that 
\begin{eqnarray}
  \sqrt{Q}dz \in i\R_-, \qquad \text{on $\Gamma_2\cap\C_+$} \\
  \sqrt{Q}dz \in i\R_+, \qquad \text{on $\Gamma_2\cap\C_-$}
\end{eqnarray} 
with the orientation from $\xi$ to $\overline{\xi}$.
From \eqref{e-hdef}, $i(2h-W')=\sqrt{Q}$.
Thus we have 
\begin{eqnarray}  
  (2h-W')dz <0, \quad z\in\Gamma_2\cap\C_+, \\
  (2h-W')dz >0, \quad z\in\Gamma_2\cap\C_-.
\end{eqnarray}
This proves condition (d).
The reader will observe that the remaining conditions and formulae
in Proposition \ref{prop-h} have been proved en route 
in this section, and this completes the proof of the Proposition.

\begin{figure}[h]
\begin{center}
\mbox{\psfig{file=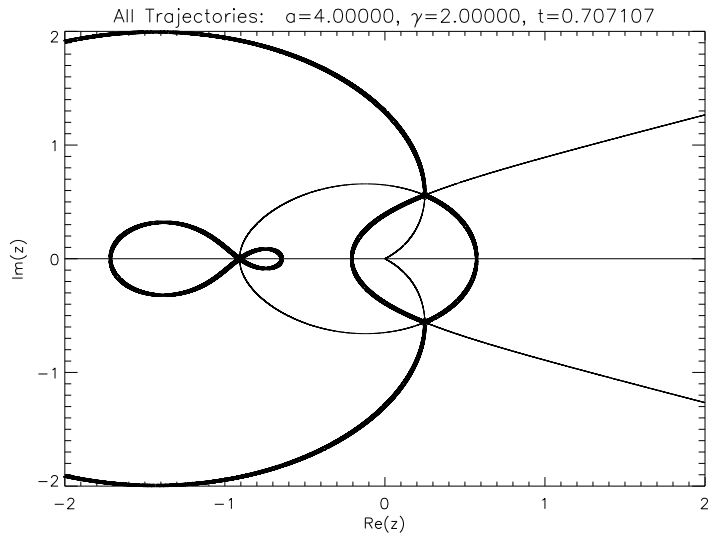,width=4 in}}
\end{center}
\caption{The trajectories and orthogonal trajectories of the quadratic
differential $Q(z)\,(dz)^2$ emerging from the points $z=\xi$,
$z=\overline{\xi}$, and $z=-z_0$.  The parameter values are $a=4$,
$\gamma=2$, and $t=1/\sqrt{2}$.  The trajectories where
$Q(z)\,(dz)^2>0$ are shown with thick curves and the orthogonal
trajectories where $Q(z)\,(dz)^2<0$ are shown with thin curves.}
\label{fig:all_trajectories}
\end{figure}

These facts can be illustrated by numerical computations of
trajectories and orthogonal trajectories associated with the quadratic
differential $Q(z)\,(dz)^2$.  The (orthogonal) trajectories can be
obtained simply with a Runge-Kutta scheme.  For (orthogonal)
trajectories that emerge from zeros or poles of $Q(z)$, which amount
to singularities of the vector field in the complex plane, the only
additional difficulty is to determine the initial direction, which is
not unique.  But the possible initial directions are easily determined
by the sort of local analysis that has already been presented above.
An example of the results of such a calculation is presented in
Figure~\ref{fig:all_trajectories}.

\begin{figure}[h]
\begin{center}
\mbox{\psfig{file=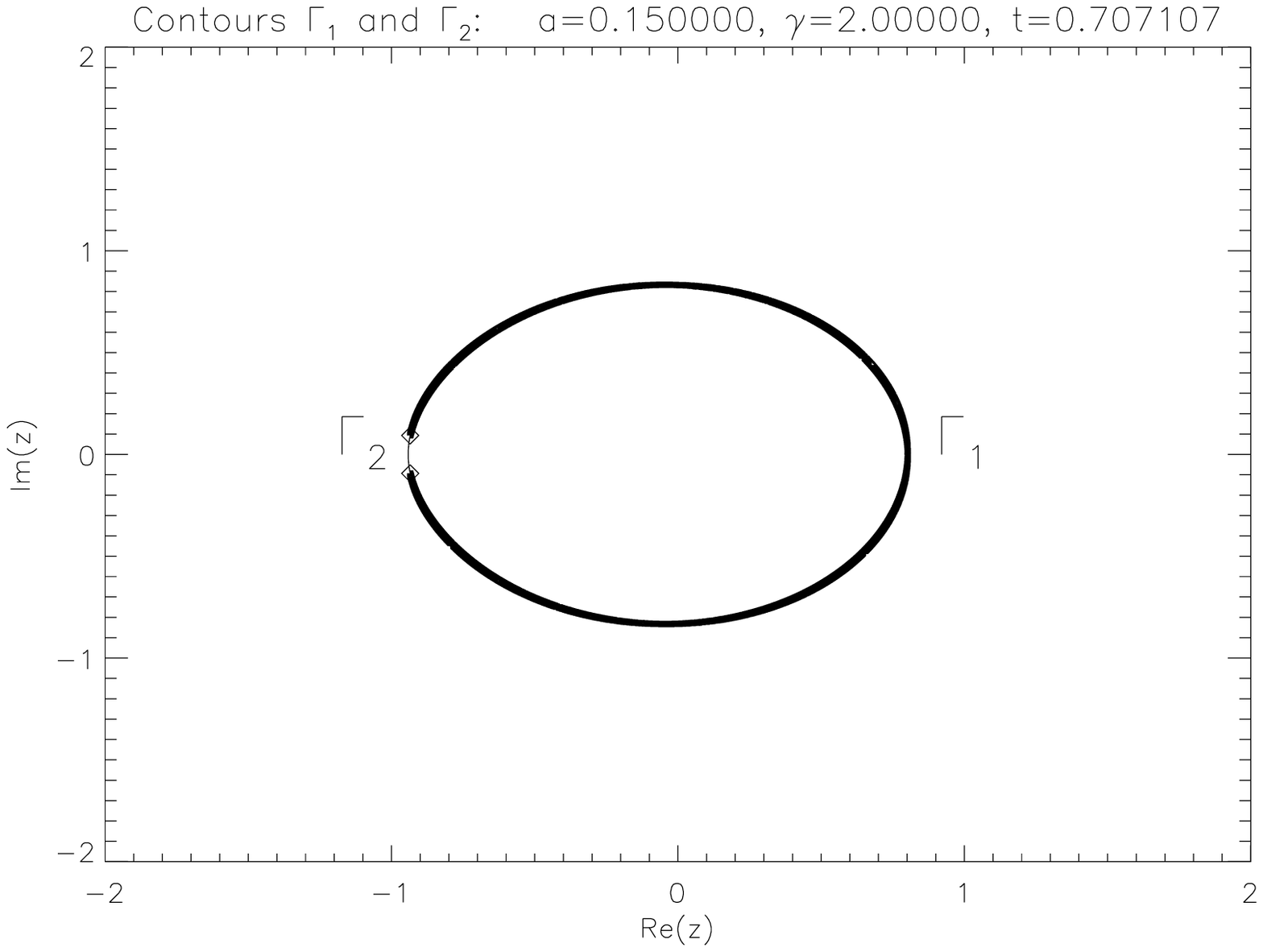,width=2 in}
\psfig{file=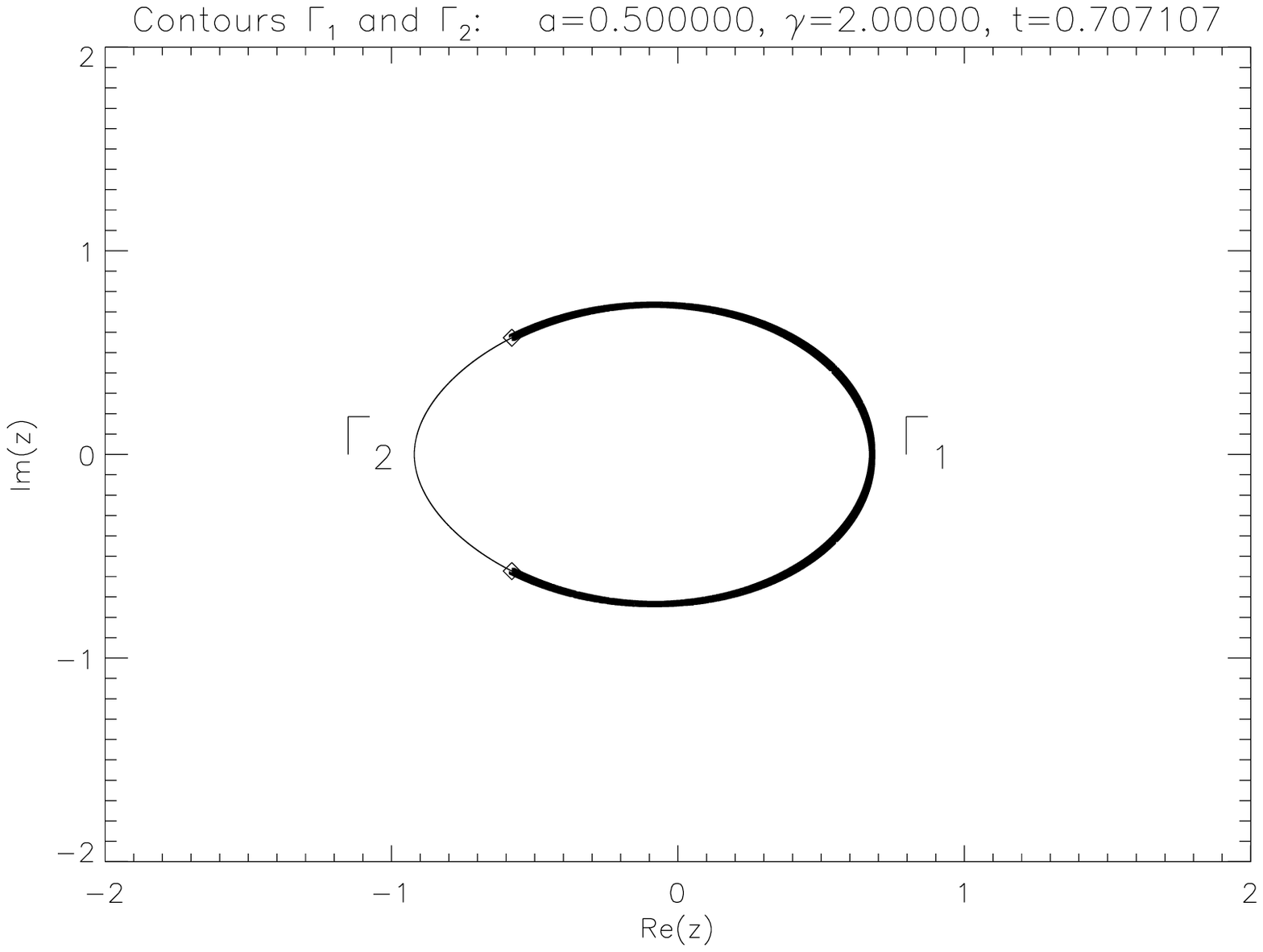,width=2 in}
\psfig{file=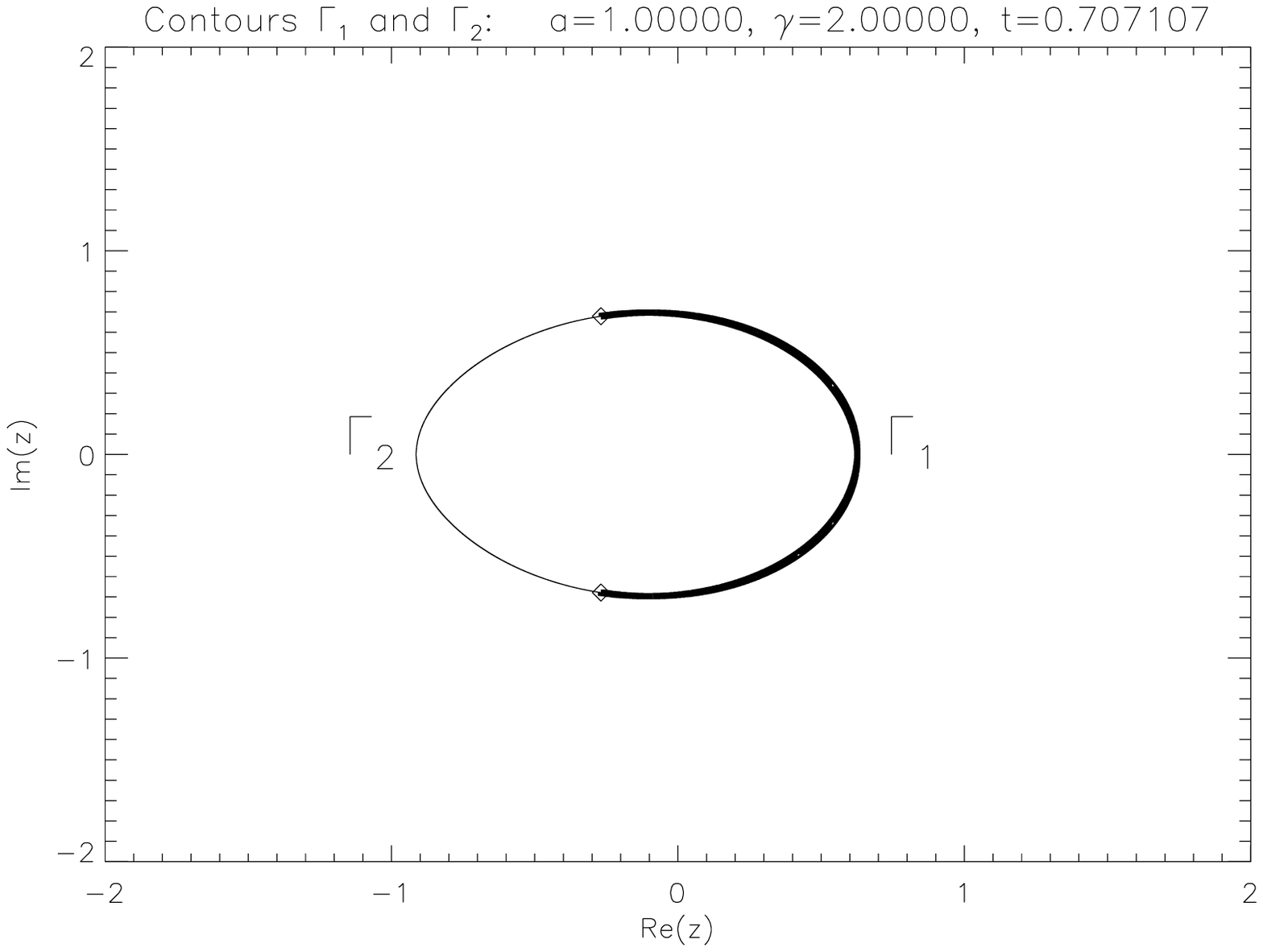,width=2 in}}\\
\mbox{\psfig{file=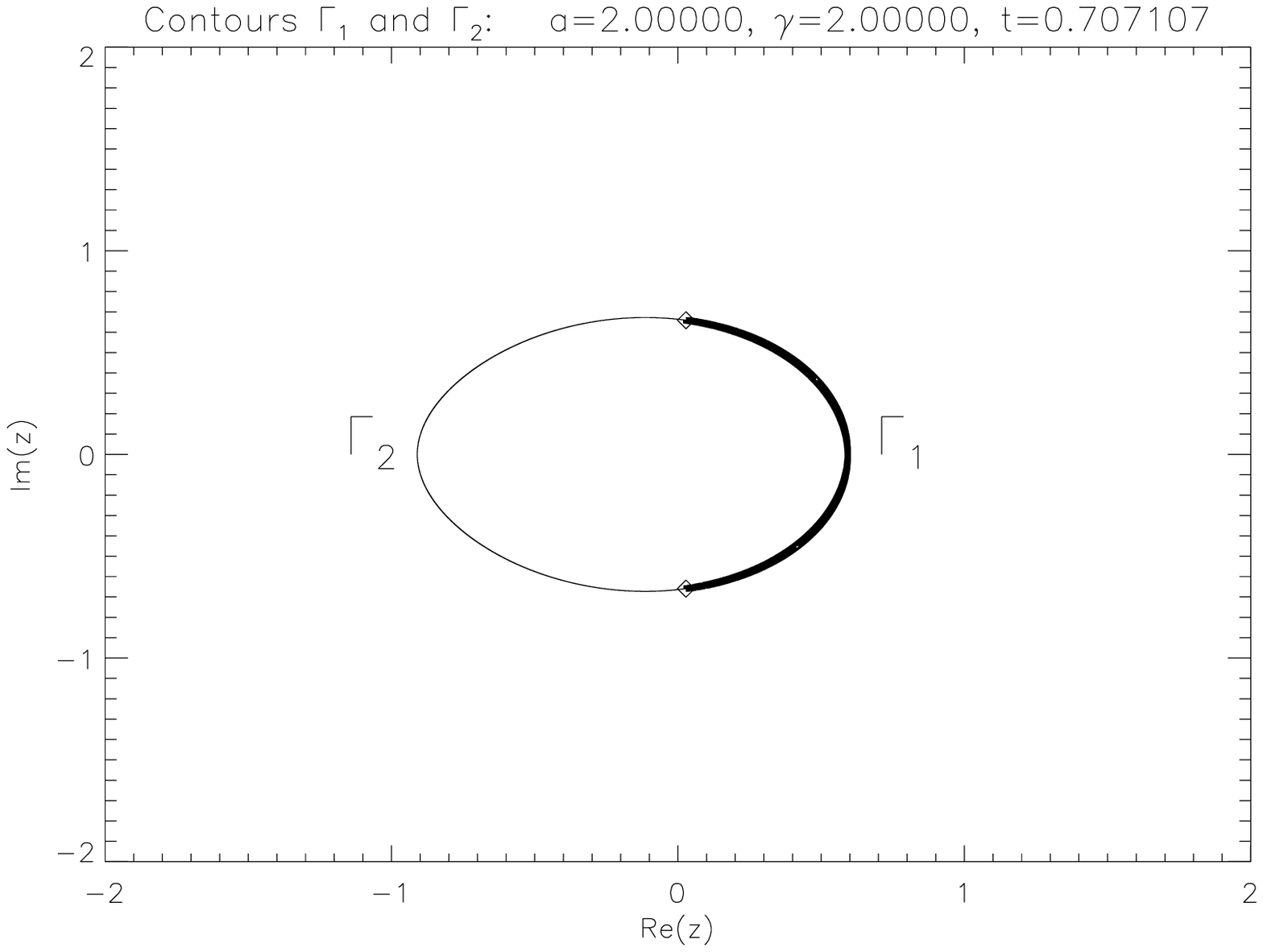,width=2 in}
\psfig{file=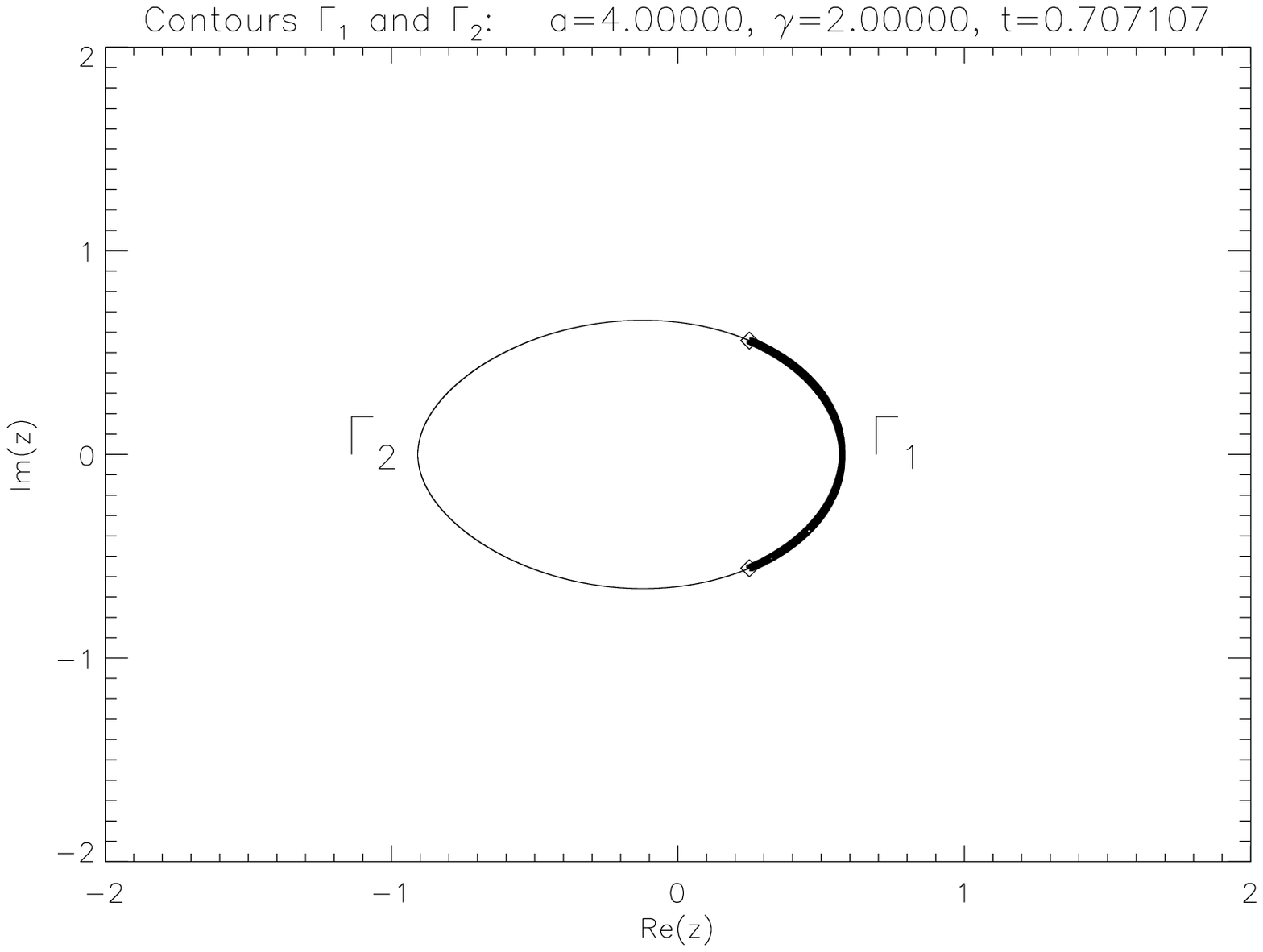,width=2 in}
\psfig{file=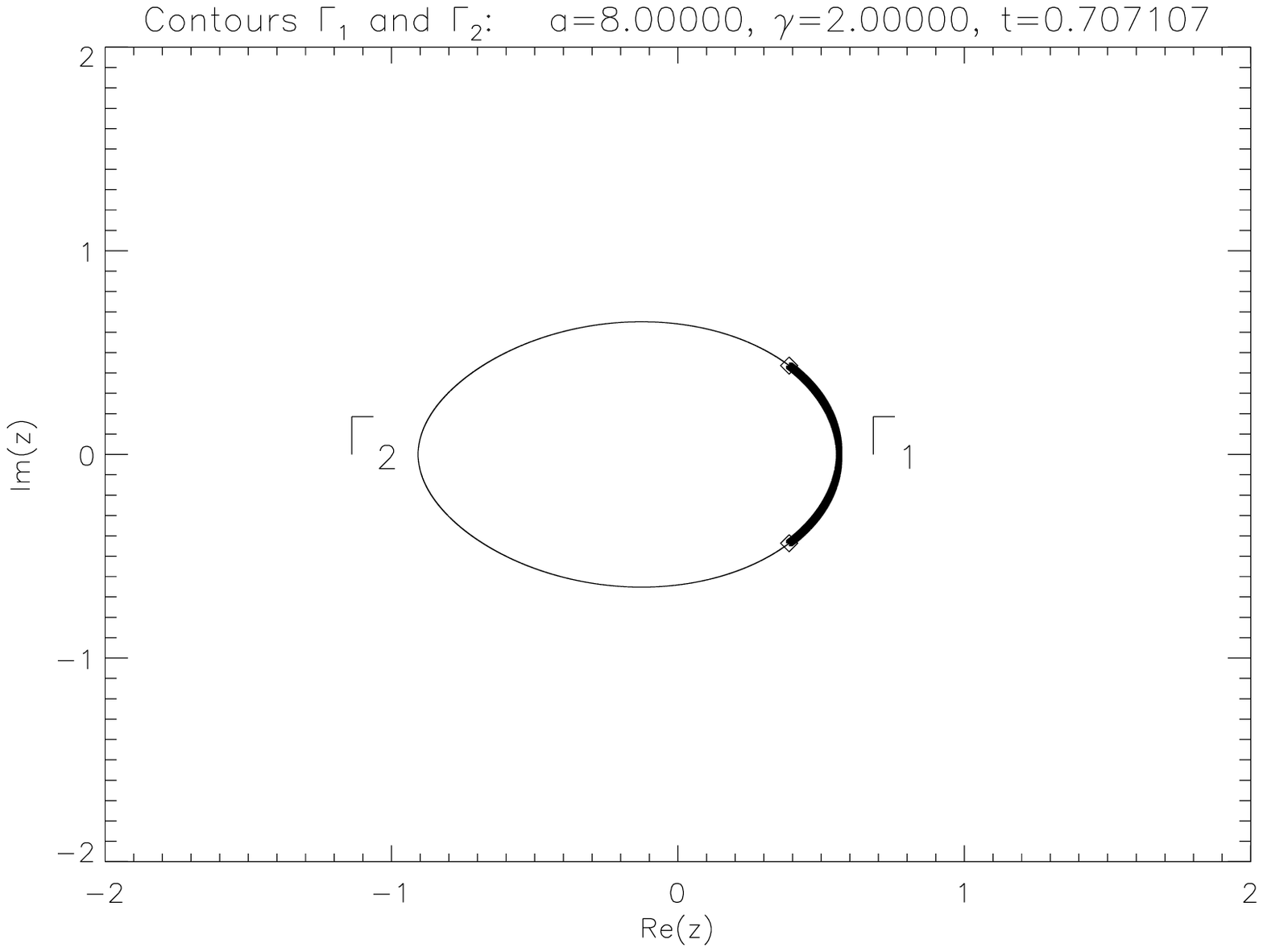,width=2 in}}
\end{center}
\caption{The dependence of the contours $\Gamma_1$ and $\Gamma_2$ on the
parameter $a$.  As $a$ tends to $a_0$ from above, the endpoints $\xi$ and
$\overline{\xi}$ coalesce on the negative real axis.}
\label{fig:contours_aseq}
\end{figure}

A computer program for generating numerical approximations to the
contours $\Gamma_1$ and $\Gamma_2$ is of course useful 
because it allows one to explore/illustrate the dependence 
of the contours on the parameters $a$, $\gamma$, and $t$.  As an example,
Figure~\ref{fig:contours_aseq} 
illustrates the deformation of the
contours as $a$ is varied while $\gamma$ and $t$ are held fixed.
Here, we see clearly what happens as $a$ is decreased to $a_0$, its
minimum value which depends on $\gamma$ and $t$.  Namely, the two
endpoints converge to a common point on the negative real axis, and at
the same time, the contour $\Gamma_1$ closes without collapsing, while
$\Gamma_2$ disappears.  Similarly, as $a$ increases without bound, the
opposite situation prevails, with the endpoints $\xi$ and
$\overline{\xi}$ coalescing on the positive real axis, while
$\Gamma_2$ closes without collapsing and $\Gamma_1$ disappears.  It is
also possible to see in these pictures that the closed curve
$\Gamma_1\cup\Gamma_2$ deforms somewhat throughout this process; the
endpoints $\xi$ and $\overline{\xi}$ do not simply slide along a fixed
closed curve as $a$ is varied.

\section{The $g$-function}\label{sec-g}

\begin{prop}\label{prop-g}
Fix $0<t<1$, $\gamma\ge 1$, $a>a_0$. 
Let $\xi$, $\Gamma_1$, $\Gamma_2$ be as in the Proposition \ref{prop-h}.
Let $\ln z$ denote the principal branch 
of logarithm, $\ln z\in \R$ for $z>0$, 
and set 
\begin{equation}\label{e-g}
  g(z):= \int_{\infty}^z \bigl(h(s)-\frac1s\bigr) ds + \ln z,
\qquad z\in\C\setminus (\overline{\Gamma_1}\cup (-\infty, p_{i}]),
\end{equation}
where $p_i>0$ is the intersection point of $\Gamma_1$ and $\R$, 
and the integral is taken over a curve from $\infty$ to $z$ 
which does not intersect $\overline{\Gamma_1}\cup (-\infty, p_{i}]$.
Let $\Phi(z)$ be as in \eqref{e-Phi} of Proposition \ref{prop-h}.
Also set 
\begin{equation}\label{e-l}
  \ell := 2g(\xi)-W(\xi).
\end{equation}
Then $g$ and $\ell$ satisfy 
the following properties:
\begin{itemize}
\item[(1)] $g(z)$ is well-defined and analytic in 
$\C\setminus (\overline{\Gamma_1}\cup (-\infty, p_{i}])$, 
and $e^{g(z)}$ is analytic in $\C\setminus\overline{\Gamma_1}$.
\item[(2)] $e^{g(z)}= z(1 + O(1/z))$ as $z\to\infty$.
\item[(3)] $e^{g_+(z)+g_-(z)-W(z)-\ell} =1$ for $z\in\Gamma_1$.
\item[(4)] $e^{g_+(z)-g_-(z)} = e^{\Psi_{1+}(z)}$ for $z\in\Gamma_1$, 
where 
$e^{\Psi_{1}(z)}=\exp\bigl\{\int_\xi^z\Phi(s)ds\bigr\}$,
$z\in \C\setminus\overline{\Gamma_1\cup\R_-}$
with the integral taken along 
any curve not intersecting $\overline{\Gamma_1\cup\R_-}$.
\item[(5)] $e^{2g(z)- W(z)-\ell} = e^{\Psi_2(z)}$ for $z\in\Gamma_2$, 
where $\Psi_2(z)= \int_\xi^z\Phi(s)ds$
with the integral taken along $\Gamma_2$.
\end{itemize}
\end{prop}

\begin{rem}\label{rem-real}
It follows from \eqref{e-hprop} that $h(z)=\overline{h(\overline{z})}$ 
for $z\in\C\setminus\overline{\Gamma_1}$, 
and hence from \eqref{e-g}, we see that 
$g(z)=\overline{g(\overline{z})}$ for 
$z\in\C\setminus(\overline{\Gamma_1}\cup (-\infty, p_i))$. 
\end{rem}

\begin{proof}
Since $h$ is analytic in $\C\setminus \overline{\Gamma_1}$ and 
continuous up to the boundary, 
we have $g$ given by \eqref{e-g} is well-defined and analytic in 
$\C\setminus \overline{\Gamma_1}\cup (-\infty, p_i])$, 
and continuous up to the boundary.
Now let $C$ be a simple closed curve, oriented counter-clockwise 
enclosing $\Gamma_1$.
Using the formula \eqref{e-hprop} for $h$, 
\begin{equation}
   \int_{C} h(z) dz = \int_{C} \frac{R(z)}{2\pi i} dz
\int_{\Gamma_1} \frac{W'(s)}{R_+(s)(s-z)} ds
= \frac1{2\pi i} \int_{\Gamma_1} \frac{W'(s)}{R_+(s)} ds   
\int_{C} \frac{R(z)}{s-z} dz.
\end{equation}
Using the residue at infinity, we then have 
\begin{equation}
  \frac1{2\pi i} \int_{C} \frac{R(z)}{s-z} dz 
= s - \frac{\xi+\overline{\xi}}2, 
\end{equation}
and 
\begin{equation}\label{e-hintGamma}
  \int_{C} h(z) dz = \int_{\Gamma_1} 
\frac{W'(s)}{R_+(s)} \bigl(s - \frac{\xi+\overline{\xi}}2\bigr) ds 
= -2\pi i
\end{equation}
from the endpoint conditions \eqref{e-endprop}.
Now from the definition of $g$, 
for $z\in (-\infty,p_i)$, $g_+(z)-g_-(z)=\int_{C} h(s)ds = -2\pi i$, 
where $g_\pm(z)= \lim_{\epsilon\downarrow 0} g(z\pm \epsilon i)$.
Therefore $e^{g(z)}$ is analytic in $\C\setminus\overline{\Gamma_1}$ 
and continuous up to the boundary.
This proves property (1).

Since $h(z) = \frac1z+O(z^{-2})$ as $z\to\infty$, 
we have $g(z) = \ln z + O(z^{-1})$ as $z\to\infty$, 
which proves property (2).

For $z\in\Gamma_1\cap\C_+$, 
\begin{equation}\label{e-gpm1}
 g_\pm(z) = g(\xi) + \int_{\xi}^z \bigl( h_\pm(s)- \frac1s \bigr) ds 
+ \ln z - \ln \xi
=  g(\xi) + \int_{\xi}^z h_\pm(s) ds,
\end{equation}
where the integral from $\xi$ to $z$ is taken along $\Gamma_1$.
Hence from $h_+(z)+h_-(z)= W'(z)$, $z\in\Gamma_1$, 
we have, using \eqref{e-l},
\begin{equation}
  g_+(z)+g_-(z) = 2g(\xi) + \int_{\xi}^z W'(s) ds
= W(z) + \ell,
\end{equation}
for $z\in\Gamma_1\cap\C_+$.
For $z\in\Gamma_1\cap\C_-$,
\begin{equation}
 g_\pm(z) = g(\overline{\xi}) 
+ \int_{\overline{\xi}}^{z} h_\pm(s) ds 
= g(\overline{\xi})
+ \int_{\overline{\xi}}^{\xi} h_\pm(s) ds 
+ \int_{\xi}^{z} h_\pm(s) ds,
\end{equation}
where the integrals are again taken along $\Gamma_1$.
But from \eqref{e-hintGamma} 
(recall that $h$ is analytic in $\C\setminus \overline{\Gamma_1}$),
$\int_{\overline{\xi}}^{\xi}  (h_+(s)-h_-(s)) ds= 2\pi i$ 
and $g(\overline{\xi})-g(\xi)= \int_{\xi}^{\overline{\xi}} h_-(s)ds$.
Hence for $z\in\Gamma_1\cap\C_-$,  
\begin{equation}\label{e-gpm2}
 g_\pm(z) = g(\xi) - \int_{\overline{\xi}}^{\xi} h_-(s)ds
+ \int_{\overline{\xi}}^{\xi} h_\pm(s)ds
+ \int_{\xi}^{z} h_\pm(s) ds,
\end{equation}
and hence 
\begin{equation}
 g_+(z)+g_-(z) = 2g(\xi) + \int_{\overline{\xi}}^{\xi}  (h_+(s)-h_-(s)) ds
 + \int_\xi^z W'(s)ds
= W(z) + \ell + 2\pi i.
\end{equation}
Therefore, 
since $e^{g_+}$, $e^{g_-}$ and $W$ are continuous for $z\in\Gamma_1$, 
we have $e^{g_+(z)+g_-(z)-W(z)-\ell} =1$ for all $z\in \Gamma_1$, 
which verifies property (3).

From \eqref{e-hPhi1}, \eqref{e-gpm1}, \eqref{e-gpm2}, 
and the above relation 
$\int_{\overline{\xi}}^{\xi}  (h_+(s)-h_-(s)) ds= 2\pi i$,
we have for $z\in\Gamma_1$, 
\begin{equation}
  e^{g_+(z)-g_-(z)} = \exp\biggl\{\int_{\xi}^z (h_+(s)-h_-(s))ds\biggr\}
= e^{\int_\xi^z \Phi(s)_+ds },
\end{equation}
and hence the property (4) follows if we prove that 
$\exp\bigl\{ \int_\xi^z \Phi(s)ds\bigr\}$, 
$z\in\C\setminus\overline{\Gamma_1\cup\R_-}$ does not depend on 
the choice of the integration path.
For this purpose, it is enough to prove that 
\begin{equation}\label{e-2piZ}
  \int_{C} \Phi(z)dz \in 2\pi i \Z
\end{equation}
for any simple closed contour $C$ which encloses $\Gamma_1$ 
and does not intersect $(-\infty, 0]$.
From \eqref{e-hPhi} and the fact that $W'(z)$ is analytic 
away from $-t^{-1}, -t, 0$, 
we have 
\begin{equation}
 \int_{C} \Phi(z)dz 
= \int_{C} (2h(z)- W'(z))dz 
= 2\int_C h(z)dz.
\end{equation}
Using \eqref{e-hprop} for $h$, the above integral is equal to 
\begin{equation}
 \int_{C} \frac{R(z)}{\pi i} dz \int_{\Gamma_1} \frac{W'(s)}{R_+(s)(s-z)} ds
= \frac1{\pi i}\int_{\Gamma_1} \frac{W'(s)}{R_+(s)} ds 
\int_C \frac{R(z)}{s-z} dz 
= 2\int_{\Gamma_1} \frac{W'(s)}{R_+(s)} 
\bigl( \frac{\xi+\overline{\xi}}2 -s \bigr) ds.
\end{equation}
But this is equal to $4\pi i$ from the endpoint condition \eqref{e-endprop},
and so \eqref{e-2piZ} is established.

For $z\in\Gamma_2\cap\C_+$, as in \eqref{e-gpm1}
\begin{equation}
  g(z) = g(\xi) + \int_\xi^z h(s)ds
\end{equation}
and 
\begin{equation}
  2g(z) - W(z) -\ell = 2\int_\xi^z h(s)ds - W(z) + W(\xi) 
= \int_\xi^z (2h(s) -W'(s)) ds.
\end{equation}
For $z\in\Gamma_2\cap\C_-$, as in \eqref{e-gpm2}, 
\begin{equation}
  g(z) = g(\overline{\xi}) + \int_{\overline{\xi}}^z h(s)ds
= g(\xi) - \int_{\overline{\xi}}^\xi h_-(s) ds 
+ \int_{\overline{\xi}}^\xi h_+(s) + 
\int_\xi^z h(s)ds
= g(\xi) + \int_\xi^z h(s)ds + 2\pi i,
\end{equation}
and hence for $z\in\Gamma_2$, 
\begin{equation}
  e^{2g(z)-W(z)-\ell} = \exp\biggl\{ \int_{\xi}^z (2h(s)-W'(s))ds \biggr\},
\end{equation}
and property (5) follows from \eqref{e-hPhi2}.
This completes the proof of Proposition \ref{prop-g}.
\end{proof}

\section{RHP analysis}\label{sec-rhp}

Set 
\begin{equation}
  \Gamma = \overline{\Gamma_1\cup\Gamma_2},
\end{equation}
oriented counter-clockwise.
It is a simple closed curve which has $0$ and $-t$ inside and 
$-t^{-1}$ outside.
Since the jump matrix 
$V_Y=\bigl( \begin{smallmatrix} 1 & z^{-k}\varphi(z) \\ 0 & 1 
\end{smallmatrix} \bigr)$ 
for $Y$ in \eqref{e-Y} is analytic in 
$\C\setminus \{0\}$, we can deform the contour $\Sigma$ for $Y$ to $\Gamma$,
as follows.
Set 
\begin{equation}
\widetilde{Y}(z)=
\begin{cases}
  Y(z) \qquad \text{for $z$ inside both $\Gamma$ and $\Sigma$, 
and for $z$ outside both $\Gamma$ and $\Sigma$,} \\
  Y(z)V_Y(z) \qquad \text{for $z$ inside $\Gamma$ and outside $\Sigma$,} \\
  Y(z)V_Y^{-1}(z) \qquad \text{for $z$ outside $\Gamma$ and inside $\Sigma$.}
\end{cases}
\end{equation}
Then $\widetilde{Y}$ is analytic in $\C\setminus\Gamma$ and continuous 
up to the boundaries, satisfies $\widetilde{Y}_+(z)=\widetilde{Y}_-V_Y(z)$ 
for $z\in \Gamma$, and $\widetilde{Y}(z)z^{-k\sigma_3}= I+O(z^{-1})$ 
as $z\to\infty$.

Now (see the Introduction) we define 
\begin{equation}
  M(z) = e^{-\frac12\ell k \sigma_3} \widetilde{Y}(z) 
e^{-(g(z)-\frac12\ell)k\sigma_3}.
\end{equation}
Then from Proposition \ref{prop-g} (1), (2), 
$M$ satisfies 
\begin{equation}\label{e-M}
\begin{cases}
  M(z) \text{ is analytic in $z\in\C\setminus \Gamma$}, \\
  M_+(z)= M_-(z) \begin{pmatrix} e^{-k(g_+(z)-g_-(z))}
& e^{k(g_+(z)+g_-(z)-W(z)-\ell)} \\ 0&e^{k(g_+(z)-g_-(z))} 
\end{pmatrix},
\qquad z\in\Gamma, \\
  M(z) = I + O(1/z), \qquad \text{as $z\to\infty$}.
\end{cases}
\end{equation}
From the Proposition \ref{prop-g} (3), (4), the jump matrix $V$ for $M$ is now 
\begin{equation}
  V(z) = \begin{pmatrix} e^{-k\Psi_{1+}(z)} & 1 \\ 0 & 
e^{k\Psi_{1+}(z)} 
\end{pmatrix},
\qquad z\in\Gamma_1, 
\end{equation}
and from the Proposition \ref{prop-g} (1), (5), we have 
\begin{equation}
  V(z) = \begin{pmatrix} 1 & e^{k\Psi_2(z)} \\ 0 & 1 \end{pmatrix},
\qquad z\in\Gamma_2.
\end{equation}

\begin{figure}[ht]
 \centerline{\epsfig{file=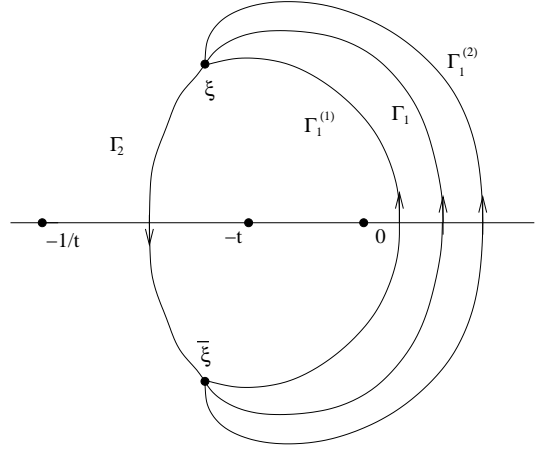, width=7cm}}
 \caption{The contours $\Gamma_1^{(1)}$ and $\Gamma_1^{(2)}$.}
\label{fig-deformed}
\end{figure}

For the jump matrix on $z\in\Gamma_1$, note that $\Psi_+(z)=-\Psi(z)_-$ and 
\begin{equation}
V(z)= 
\begin{pmatrix} e^{-k\Psi_{1+}(z)} & 1 \\ 0 & e^{-k\Psi_{1-}(z)} 
\end{pmatrix}
= \begin{pmatrix} 1&0 \\ e^{-k\Psi_{1-}(z)} & 1 \end{pmatrix}
\begin{pmatrix} 0 & 1 \\ -1 & 0 \end{pmatrix}
\begin{pmatrix} 1&0 \\ e^{-k\Psi_{1+}(z)} & 1 \end{pmatrix}. 
\end{equation}
Clearly $\Psi_{1+}$ has an analytic continuation to the $(+)$-side 
of the contour $\Gamma_1$. 
Now for $z\in\Gamma_1$, 
it is easy to see that $Re(\Psi_{1_+}(z))=0$, and hence
from the Proposition \ref{prop-h} (c) and \eqref{e-hPhi1}, 
the derivative $\frac{d}{dz} \Psi_{1+}(z)= \Phi_+(z)$ 
along the contour $\Gamma_1$ satisfies
$Im( \frac{d}{dz}\Psi_{1+}(z))<0$, $z\in\Gamma_1$.
Thus the Cauchy-Riemann condition for the analyticity
implies that $Re(\Psi_{1}(z))>0$ for $z$ on the $(+)$-side of $\Gamma_1$ 
and close to the contour.
Therefore we can take a contour $\Gamma_1^{(1)}$ with endpoints
$\xi$, $\overline{\xi}$ for which $Re(\Psi(z))>0$ 
for $z\in int(\Gamma^{(1)})$.
Similarly, $\Psi_{1-}$ has an analytic continuation to the $(-)$-side 
of $\Gamma_1$ and its real part is positive for 
$z$ on the $(+)$-side of $\Gamma_1$ close to the contour $\Gamma_1$
and we take a contour $\Gamma_1^{(2)}$ for which $Re(\Psi(z))>0$ 
on its interior.
We take the orientation of $\Gamma_1^{(j)}$, $j=1,2$ to be from 
$\overline{\xi}$ to $\xi$.
See Figure \ref{fig-deformed} for the general shape of the contours 
$\Gamma_1^{(1)}$, $\Gamma_1^{(2)}$.

Define $\widetilde{M}(z)$ to be $M(z)$ for $z$ in the region bounded by 
$\Gamma_2$ and $\Gamma_1^{(1)}$ and also in the unbounded region.
For the region bounded by $\Gamma_1$ and $\Gamma_1^{(1)}$, define 
\begin{equation}
\widetilde{M}= M
\begin{pmatrix} 1&0 \\ -e^{-k\Psi_{1}(z)} & 1 \end{pmatrix},
\end{equation}
and for the region bounded by $\Gamma_1$ and $\Gamma_1^{(2)}$, define 
\begin{equation}
\widetilde{M}= M
\begin{pmatrix} 1&0 \\ e^{-k\Psi_{1}(z)} & 1 \end{pmatrix}.
\end{equation}
Set $\Gamma'=\overline{\Gamma_1\cup\Gamma_1^{(1)}
\cup\Gamma_1^{(2)}\cup\Gamma_2}$.
Then $\widetilde{M}$ satisfies the new RHP 
\begin{equation}\label{e-wideM}
\begin{cases}
  \widetilde{M}(z) \text{ is analytic in $z\in\C\setminus \Gamma'$}, \\
  \widetilde{M}_+(z)= \widetilde{M}_-(z) \widetilde{V}(z)
\qquad z\in\Gamma', \\
  \widetilde{M}(z) = I + O(1/z), \qquad \text{as $z\to\infty$},
\end{cases}
\end{equation}
where the jump matrix $\widetilde{V}$ is 
\begin{equation}
\widetilde{V}(z)=
\begin{cases}
\begin{pmatrix} 0&1\\-1&0 \end{pmatrix}, 
\qquad &z\in\Gamma_1, \\
\begin{pmatrix} 1&0 \\ e^{-k\Psi_{1}(z)} & 1 \end{pmatrix}, 
\qquad &z\in\Gamma_1^{(1)}\cup\Gamma_1^{(2)}, \\
\begin{pmatrix} 1& e^{k\Psi_{2}(z)} \\ 0 & 1 \end{pmatrix},
\qquad &z\in\Gamma_2.
\end{cases}
\end{equation}

Now we take the limit $k\to\infty$ with $\gamma\ge 1$ in a compact set 
and $a_0<a\le a^*$, for some $a^*$.
From the signature of $Re(\Psi_1(z))$ on $\Gamma_1^{(j)}$, $j=1,2$, 
we see that the jump matrix 
$\widetilde{V} \to I$.
For $z\in\Gamma_2$, we have $\Psi_2(z)<0$
from the Proposition \ref{prop-g} and \ref{prop-h}.
Indeed $\Psi_2(z)$ is decreasing as $z$ follows from $\xi$ to $-z_0$ 
along $\Gamma_2$, and then increasing as $z$ follows from $-z_0$ to 
$\overline{\xi}$ along $\Gamma_2$.
On the other hand, 
$\Psi_2(\xi)=0$ and 
$\Psi_2(\overline{\xi})=\int_{\xi}^{\overline{\xi}}\Phi(s)ds =0$,
as $\Phi(s)ds=\overline{\Phi(s)ds}=-\Phi(\overline{s})d\overline{s}$,
and the negativity of $\Psi_2(z)$ on $\Gamma_2$ follows.
Hence as $k\to\infty$, $\widetilde{V} \to I$ on $\Gamma_2$.
Therefore we have $\widetilde{V} \to V^\infty$ where
\begin{equation}
 V^\infty= 
\begin{cases}
\begin{pmatrix} 0&1\\-1&0 \end{pmatrix}, \qquad &z\in \Gamma_1,\\
I \qquad &z\in\Gamma_1^{(1)}\cup\Gamma_1^{(2)}\cup\Gamma_2.
\end{cases}
\end{equation}
Let $M^\infty$ be the solution to the RHP with the jump matrix $V^\infty$
and normalized at infinity.
The solution is given by
\begin{equation}
  M^\infty(z) = \begin{pmatrix}
\frac{\beta+\beta^{-1}}2 & \frac{\beta-\beta^{-1}}{2i} \\
- \frac{\beta-\beta^{-1}}{2i} & \frac{\beta+\beta^{-1}}2 
\end{pmatrix},
\end{equation}
where
\begin{equation}
  \beta(z) = \biggl( \frac{z-\xi}{z-\overline{\xi}} \biggr)^{1/4}
\end{equation}
which is defined to be analytic $\C\setminus\Gamma_1$ 
and satisfies $\beta(z)\sim 1$ as $z\to\infty$.
We expect that $\widetilde{M} \sim M^\infty$ as $k\to\infty$, and hence 
by tracking the algebraic transformations $Y\to M\to \widetilde{M}$, 
we expect that
\begin{equation}\label{e-YMinfty}
  Y_{21}(0)e^{-k(g(0)-\ell)} \sim M^\infty_{21}(0), 
\qquad k\to\infty,
\end{equation}
for $\gamma\ge 1$ in a compact set and $a_0<a\le a^*$.

Indeed we have:
\begin{prop}
Let $1\le \gamma\le \gamma_1$ for any fixed $\gamma_1\ge 1$.
There are $L_1, \delta_1>0$ such that 
for 
\begin{equation}
  a_0 + \frac{L_1}{k^{2/3}} \le a\le (1+\delta_1)a_0,
\end{equation}
we have 
\begin{equation}\label{e-YMinfty1}
  Y_{21}(0)e^{-k(g(0)-\ell)} = M^\infty_{21}(0) 
\bigl( 1+ O(\frac1{k|a-a_0|})\bigr),
\end{equation}
for sufficiently large $k$.
\end{prop}

The convergence $\widetilde{V} \to V^\infty$ is not uniform on
$\Gamma'$ and this considerably complicates the analysis.  As in
\cite{DKMVZ3} in order to obtain the above error bound, we need to
introduce local parametrices for the solution of the RHP near each of
the endpoints $\xi$ and $\overline{\xi}$.  As in \cite{BDJ}, a
suitable local parametrix near each endpoint can be obtained in terms
of Airy functions.  Also since $a$ is not fixed, but is allowed to
approach $a_0$, we need to vary the magnitude of the parametrix
according to the size of $a-a_0$.  
A similar situation
arises in Lemma 6.2 (ii) of \cite{BDJ}.  The proof of the above
Proposition is parallel to that of Lemma 6.2 (ii), \cite{BDJ}
($\gamma$, $q$ in \cite{BDJ} play the same role as $a$, $k$ in this
paper, respectively), and we do not repeat the argument here.

\bigskip

\begin{lem}\label{lem-real}
We have $\Delta:=g(0+i0)-\ell\in \R$. 
In particular, $e^{g(0)-\ell}=e^{g(0+i0)-\ell}>0$.
\end{lem}

\begin{proof}
For $x\in\R\setminus\{p_i\}$, by Remark \ref{rem-real}, 
$g(x+i0)=\overline{g(x-i0)}$, and hence by Proposition \ref{prop-g} (1), 
$e^{g(x+i0)}=e^{g(x-i0)}=e^{\overline{g(x+i0)}}$,
and so $e^{g(x)}=e^{g(x+i0)}=e^{g(x-i0)}$ is real. 
In particular, it follows that $Im(g(x+i0))\in \Z \pi$, 
and hence by continuity, $Im(g(x+i0))$ is constant for $x<p_i$.
From \eqref{e-l} and the proof of Proposition \ref{prop-g} (3), 
$\ell=g_+(\xi)+g_-(\xi)-W(\xi)= g_+(p_i+i0)+g_-(p_i)-W(p_i)$. 
Hence $e^{g(0)-\ell}=e^{g(0+i0)-g_+(p_i+i0)}e^{-g_-(p_i)+W(p_i)}$.
But by the above, $Im(g(0+i0))=Im(g_+(p_i+i0))$.
Hence $g(0+i0)-g_+(p_i+0i)\in\R$.
Clearly $g_-(p_i)$ and $W(p_i)$ are also real, and this proves the lemma.
\end{proof}

If we set $\xi=|\xi|e^{i\theta_c}$, $0<\theta_c<\pi$, then we can check 
$\beta(0)= e^{i\theta_c/2}$ and 
$M^\infty_{21}(0)= -\frac1{2i}(\beta(0)-\beta(0)^{-1}) 
= -\sin \frac{\theta_c}2$.
Hence the above proposition yields that 
\begin{equation}\label{e-Yest}
  -Y_{21}(0) = e^{k\Delta} 
\sin\frac{\theta_c}2 \bigl(1+O(\frac1{k|a-a_0|})\bigr).
\end{equation}
Note that $-Y_{21}(0)$ is indeed real and positive from Lemma \ref{lem-GandY}.
This is consistent with Lemma \ref{lem-real}.

Now we compute $e^{g(0)-\ell}$.
Let 
\begin{equation}
  \alpha = \frac{\xi+\overline{\xi}}2 = |\xi|\cos\theta_c.
\end{equation}
From the formula \eqref{e-hprop} for $h$, one can check directly that
an anti-derivative of $2(h(z)-\frac1z)$ is
\begin{equation}\label{e-hantider}
\begin{split}
&-\gamma a \log(z+t^{-1})-a\log(z+t) + (a-1)\log(z) \\
&+ \biggl( -\frac{\gamma a}{x} -\frac{a}{y} + \frac{a+1}{r} \biggr) R(z)
- \biggl( -\frac{\gamma a}{x}(\alpha+t^{-1}) - \frac{a}{y}(\alpha+t)
+\frac{a+1}{r}\alpha \biggr) \log(z-\alpha+R(z)) \\
&- \gamma a \log\biggl(\frac{z+R(z)+t^{-1}-x}{z+R(z)+t^{-1}+x} \biggr)
- a \log\biggl(\frac{z+R(z)+t-y}{z+R(z)+t+y} \biggr) 
+ (a+1) \log\biggl(\frac{z+R(z)-r}{z+R(z)+r} \biggr)  
\end{split}
\end{equation}
where the logarithms are taken to be the analytic 
in $\C\setminus(-\infty, 0]$ and 
$\log z= \log|z|$ for $z>0$,
and $r=-R(0)$, $x=-R(-t^{-1})$, $y=-R(-t)$ as in \eqref{e-rxyxi}.
It is straightforward, but tedious, to check that \eqref{e-hantider} is 
analytic in $\C\setminus(\overline{\Gamma_1}\cup(-\infty, p_i])$ 
as in the definition of $g$ in \eqref{e-g}.
Using the endpoint conditions \eqref{e-endprop}, or\eqref{e-contemp}, 
\eqref{e-hantider} is equal to 
\begin{equation}
\begin{split}
&-\gamma a \log(z+t^{-1})-a\log(z+t) + (a-1)\log(z) 
+ (1+\gamma a)\log(z-\alpha+R(z)) \\
&- \gamma a \log\biggl(\frac{z+R(z)+t^{-1}-x}{z+R(z)+t^{-1}+x} \biggr)
- a \log\biggl(\frac{z+R(z)+t-y}{z+R(z)+t+y} \biggr)
+ (a+1) \log\biggl(\frac{z+R(z)-r}{z+R(z)+r} \biggr).
\end{split}
\end{equation}
Evaluating the asymptotics as $z\to\infty$, we see that 
\begin{equation}\label{e-gexpl}
\begin{split}
  2g(z) 
= & -\gamma a \log(z+t^{-1})-a\log(z+t) + (a+1)\log(z) \\
&+ (1+\gamma a)\log\bigl((z-\alpha+R(z))/2\bigr) 
- \gamma a \log\biggl(\frac{z+R(z)+t^{-1}-x}{z+R(z)+t^{-1}+x} \biggr) \\
& - a \log\biggl(\frac{z+R(z)+t-y}{z+R(z)+t+y} \biggr)
+ (a+1) \log\biggl(\frac{z+R(z)-r}{z+R(z)+r} \biggr),
\end{split}
\end{equation}
and 
\begin{equation}\label{e-g0expl}
\begin{split}
  2Re(g(0+i0))
= & (\gamma-1) a \log t
+ (1+\gamma a)\log((\alpha+r)/2) \\
&- \gamma a \log\biggl|\frac{-r+t^{-1}-x}{-r+t^{-1}+x} \biggr|
- a \log\biggl|\frac{-r+t-y}{-r+t+y} \biggr|
+ (a+1) \log\biggl(\frac{2r}{1+\alpha/r} \biggr).
\end{split}
\end{equation}
Also from \eqref{e-gexpl} and \eqref{e-l}, we have 
\begin{equation}\label{e-lexpl}
\begin{split}
  \ell = &\gamma a \log t + (1+\gamma a)\log(\xi-\alpha)/2 \\
& - \gamma a \log\biggl(\frac{\xi+t^{-1}-x}{\xi+t^{-1}+x} \biggr)
- a \log\biggl(\frac{\xi+t-y}{\xi+t+y} \biggr)
+ (a+1) \log\biggl(\frac{\xi-r}{\xi+r} \biggr).
\end{split}
\end{equation}
By \eqref{e-xiextra}, 
\begin{equation}
  \cos\theta_c= \frac{x^2-r^2-t^{-2}}{2rt^{-1}}
= \frac{y^2-r^2-t^2}{2rt}.
\end{equation}
Thus using \eqref{e-g0expl} and \eqref{e-lexpl},
we can express $\Delta$ in terms of $r, x, y$.
After some algebra, we find
\begin{equation}\label{e-Delta}
\begin{split}
  \Delta =
&-\gamma a\log t +(2+a+\gamma a)\log 2 + (1+a)\log r 
-\frac12 \log(r+t^{-1}-x) \\
&-\frac12(1+2\gamma a)\log(r+t^{-1}+x)
-\frac12 \log(r+t-y) -\frac12(1+2a)\log(r+t+y).
\end{split}
\end{equation}

\medskip

To emphasize the dependence on $a$, 
we write $\Delta=\Delta(a)$, etc.

\begin{lem}\label{lem-Delta}
Fix $1\le \gamma_2<\infty$.
Then there exists $\delta_2>0$ such that for $a_0\le a\le (1+\delta_2)a_0$ 
and $1\le \gamma\le \gamma_2$,
we have 
\begin{equation}\label{e-Delest}
  \Delta(a) = -c_2(a-a_0)^2 \bigl(1 +O(|a-a_0|\bigr)
\end{equation}
where 
\begin{equation}
  c_2=\frac{t^2(t+t\gamma+2\sqrt\gamma)^3\sqrt\gamma}
{4(1+t\sqrt\gamma)^2(t+\sqrt\gamma)^2},
\end{equation}
and the order term $O(|a-a_0|)$ is uniform for 
$a$ and $\gamma$ as above.
\end{lem}

For the proof, we need the following lemma, 
which considers the case when $a=a_0$.
This case is specifically excluded from Lemma \ref{lem-H1}.

\begin{lem}\label{lem-H2}
For fixed $0<t<1$, $\gamma\ge 1$,
and $a=a_0=\frac{1-t^2}{t((\gamma+1)t+2\sqrt\gamma)}$, 
there is a unique solution $r$ to \eqref{e-Heq} satisfying
$r_1< r <r_2$, given by 
\begin{equation}\label{e-rata0}
  r=r(a_0)=r_0:=\frac{1+t\sqrt\gamma}{t+\sqrt\gamma}.
\end{equation}
The function $a\mapsto r(a)$ is smooth for all $a\ge a_0$ and 
\begin{equation}\label{e-rpata0}
  r'(a_0)=  -\frac{3(\gamma-1)(t\gamma+t+2\sqrt\gamma)^2t^2}
{4(1+t\sqrt\gamma)(t+\sqrt\gamma)^3}.
\end{equation}
\end{lem}

\begin{proof}
When $a=a_0$, we have $r_1= tr_0^2$ and $r_2=\frac1{t}r_0^2$.
Since $t<r_0<\frac1t$, $r_0$ satisfies $r_1<r_0<r_2$.
It is then a direct calculation to check that $r_0$ satisfies the 
equation \eqref{e-Heq}.
Now we want to show the uniqueness of the solution.
Let $H$ be as in \eqref{e-H} in the proof of Lemma \ref{lem-H1}.
Now when $a=a_0$, the value $r_*$ of \eqref{e-minr0} is $r_0$.
Thus we have 
\begin{equation}
 H(r)+\frac{r}{2}H'(r) \ge H(r_0)+\frac{r_0}{2}H'(r_0)=0, 
\end{equation}
for $r_1<r<r_2$, and the inequality is strict for $r\neq r_0$.
Hence if there is a zero $r_c\neq r_0$, it should satisfy 
$H'(r_c)>0$.
On the other hand, by direct calculation, we have 
\begin{equation}
  H'(r_0)=H''(r_0)=0, 
\qquad H^{(3)}(r_0) = 
\frac{24t(t+\sqrt\gamma)^5}{\sqrt\gamma(1-t^2)^2(1+t\sqrt\gamma)^3} >0,
\end{equation}
and hence $H$ is also increasing at the zero $r=r_0$.
As in Lemma \ref{lem-H1}, there is no other zero $r_c\neq r_0$ 
in $(r_1, r_2)$.

Now consider $H=H(a, r)$.
By direct calculations, we find 
\begin{equation}
  H(a_0, r_0)=H_r(a_0, r_0)=H_a(a_0, r_0)= H_{rr}(a_0,r_0)=0
\end{equation}
and 
\begin{eqnarray}
  H_{ra}(a_0,r_0)&=&\frac{4(t+\sqrt\gamma)(t\gamma+t+2\sqrt\gamma)^2t^2}
{(1+t\sqrt\gamma)^3(1-t^2)}\neq 0, \\
H_{aa}(a_0,r_0)&=& \frac{6(\gamma-1)(t\gamma+t+2\sqrt\gamma)^4t^4}
{(1+t\sqrt\gamma)^4(t+\sqrt\gamma)^2(1-t^2)}.
\end{eqnarray}
Hence near $(a_0, r_0)$, the power series of $H(a,r)$ has the form 
\begin{equation}
  H(a,r)= H_{ra}(a_0,r_0)(r-r_0)(a-a_0) + \frac12H_{aa}(a_0,r_0)(a-a_a)^2
+ \frac16H_{rrr}(a_0,r_0)(r-r_0)^3+\cdots.
\end{equation}
Motivated by this expansion, we set 
\begin{equation}
   \eta:= \frac{r-r_0}{a-a_0},
\end{equation}
and substitute $r=r_0+\eta(a-a_0)$ in $H$, and define 
\begin{equation}
  F(a, \eta) := \frac{H(a, r_0+\eta(a-a_0))}{(a-a_0)^2}.
\end{equation}
Setting 
\begin{equation}
  \eta_0:= -\frac{H_{aa}(a_0,r_0)}{2H_{ra}(a_0, r_0)}
= -\frac{3(\gamma-1)(t\gamma+t+2\sqrt\gamma)^2t^2}
{4(1+t\sqrt\gamma)(t+\sqrt\gamma)^3},
\end{equation}
a direct calculation shows that 
$F(a, \eta)$ is a smooth function near $(a_0, \eta_0)$ and 
\begin{equation}
  F(a_0, \eta_0)=0, \qquad F_{\eta}(a_0, \eta_0)=H_{ra}(a_0, r_0)\neq 0.
\end{equation}
Therefore by the implicit function theorem, there is the smooth function 
$\eta=\eta(a)$, $a_0\le a<a_0(1+\delta)$ for some $\delta>0$ such that 
$F(a,\eta(a))=0$.
Then $H(a, r_0+\eta(a)(a-a_0))=0$, and by the uniqueness of the solution, 
$r=r_0+\eta(a)(a-a_0)$ is smooth in $a$ for $a_0\le a<a_0(1+\delta)$.
But for $a>a_0$, $H_a(a, r(a))>0$, and the smoothness of $r=r(a)$ 
is elementary.
Hence $r(a)$ is a smooth function for $a\ge a_0$.
\end{proof}

Let 
\begin{equation}\label{e-xyata0}
  x_0:= x(a_0)=\frac{\sqrt\gamma(1-t^2)}{t(t+\sqrt\gamma)}, 
\qquad y_0:=y(a_0)=\frac{1-t^2}{t+\sqrt\gamma},
\end{equation}
which are obtained by setting $a=a_0$, $r=r(a_0)=r_0$ in 
\eqref{e-con1}, \eqref{e-con2}.
Then 
\begin{equation}
  \frac{x_0^2-r_0^2-t^{-2}}{2r_0t^{-1}}=-1,
\end{equation}
and hence from \eqref{e-cosxi}, the point $\xi$ is on the negative
real line.  Thus when $a=a_0$,  the two endpoints $\xi$ and
$\overline{\xi}$ collapse to the point $-r_0$ on the real line.
This is an extreme case of the deformation illustrated in 
Figure~\ref{fig:contours_aseq}.

\begin{proof}[Proof of Lemma \ref{lem-Delta}]
When $a=a_0$, we have $r=r_0$, $x=x_0$, $y=y_0$, and we can direct check 
from \eqref{e-Delta} that $\Delta(a_0)=0$.
We have 
\begin{equation}
\begin{split}
  \Delta'(a)= 
&-\gamma \log t +(1+\gamma )\log 2 + \log r 
- \gamma\log(r+t^{-1}+x) -\log(r+t+y) \\
&+ (1+a)\frac{r'}{r} -\frac12 \frac{r'-x'}{r+t^{-1}-x} 
-\frac12(1+2\gamma a)\frac{r'+x'}{r+t^{-1}+x} 
-\frac12 \frac{r'-y'}{r+t-y} -\frac12(1+2a)\frac{r'+y'}{r+t+y}.
\end{split}
\end{equation}
At $a=a_0$, from from \eqref{e-rata0}, \eqref{e-xyata0}, we have 
\begin{equation}\label{e-Deltapptemp1}
  -\gamma \log t +(1+\gamma )\log 2 + \log r_0- \gamma\log(r_0+t^{-1}+x_0) 
-\log(r_0+t+y_0) =0
\end{equation}
and hence again from \eqref{e-rata0}, \eqref{e-xyata0}, after some algebra, 
\begin{equation}
 \Delta'(a_0)= -\frac{(\gamma-1)(1-t^2)}
{(t+t\gamma+2\sqrt\gamma)(1+t\sqrt\gamma)}
\bigl[(1+t\sqrt\gamma)r'(a_0)+t\sqrt\gamma x'(a_0)-y'(a_0)\bigr].
\end{equation}
Now from the relation \eqref{e-con3} between $r,x,y$, we have 
\begin{equation}\label{e-rxyp}
  rr'= \frac1{1-t^2}(yy'-t^2xx').
\end{equation}
This implies by \eqref{e-rata0}, \eqref{e-xyata0},
\begin{equation}\label{e-rxypata0}
  (1+t\sqrt\gamma)r'(a_0)=y'(a_0)-t\sqrt\gamma x'(a_0), 
\end{equation}
and hence 
\begin{equation}
  \Delta'(a_0)=0.
\end{equation}

Now we compute $\Delta''(a_0)$.
We have 
\begin{equation}\label{e-Deltapptemp3}
\begin{split}
  \Delta''(a)= 
&2\frac{r'}{r} +(1+a) \biggl( \frac{r''}{r}-\bigl(\frac{r'}{r}\bigr)^2\biggr)
 -\frac12 \biggl( \frac{r''-x''}{r+t^{-1}-x} 
-\bigl(\frac{r'-x'}{r+t^{-1}-x}\bigr)^2 \biggr) \\
&-2\gamma \frac{r'+x'}{r+t^{-1}+x}
-\frac12(1+2\gamma a)\biggl( \frac{r''+x''}{r+t^{-1}+x} 
- \bigl( \frac{r'+x'}{r+t^{-1}+x} \bigr)^2 \biggr) \\
& -\frac12 \biggl( \frac{r''-y''}{r+t-y} 
-\bigl( \frac{r'-y'}{r+t-y} \bigr)^2 \biggr)
-2 \frac{r'+y'}{r+t+y}
-\frac12(1+2a)\biggl( \frac{r''+y''}{r+t+y}
- \bigl( \frac{r'+y'}{r+t+y} \bigr)^2 \biggr).
\end{split}
\end{equation}
First consider the terms with double derivatives. At $a=a_0$, 
\begin{equation}\label{e-Deltapptemp2}
\begin{split}
 &\biggl( (1+a) \frac{r''}{r} -\frac12 \frac{r''-x''}{r+t^{-1}-x} 
- \frac12(1+2\gamma a) \frac{r''+x''}{r+t^{-1}+x}  
 -\frac12 \frac{r''-y''}{r+t-y}  
-\frac12(1+2a) \frac{r''+y''}{r+t+y} \biggr)(a_0) \\
&= -\frac{(\gamma-1)(1-t^2)}
{4(t+t\gamma+2\sqrt\gamma)(1+t\sqrt\gamma)}
\bigl[(1+t\sqrt\gamma)r''(a_0)+t\sqrt\gamma x''(a_0)-y''(a_0)\bigr]. 
\end{split}
\end{equation}
From \eqref{e-rxyp}, 
\begin{equation}
  rr'' - \frac{1}{1-t^2} (yy''-t^2xx'')
= -(r')^2 +\frac1{1-t^2} \bigl( (y')^2-t^2(x')^2\bigr),
\end{equation}
and hence the right-hand side of \eqref{e-Deltapptemp2} is equal to 
\begin{equation}\label{e-Deltapptemp4}
 -\frac{(\gamma-1)(1-t^2)(t+\sqrt\gamma)}
{4(t+t\gamma+2\sqrt\gamma)(1+t\sqrt\gamma)}
\bigl[ -(r'(a_0))^2 
+\frac1{1-t^2} \bigl( (y'(a_0))^2-t^2(x'(a_0))^2\bigr) \bigr]. 
\end{equation}
Thus $\Delta''(a_0)$ is given by \eqref{e-Deltapptemp3} at $a=a_0$ 
where the terms with 
double derivatives are replaced by \eqref{e-Deltapptemp4},
which involves only the first derivatives of $r, x, y$ at $a_0$. 
From \eqref{e-rpata0} and \eqref{e-con1}, \eqref{e-con2}, 
we have
\begin{eqnarray}
 r'(a_0) &=&  -\frac{3(\gamma-1)(t\gamma+t+2\sqrt\gamma)^2t^2}
{4(1+t\sqrt\gamma)(t+\sqrt\gamma)^3} \\
 x'(a_0)  &=& \frac{(t+4\sqrt\gamma +3t\gamma)(t+t\gamma+2\sqrt\gamma)^2t}
{4(1+t\sqrt\gamma)(t+\sqrt\gamma)^3} \\
 y'(a_0)  &=& \frac{(4t\sqrt\gamma +\gamma +3)(t+t\gamma+2\sqrt\gamma)^2t^2}
{4(1+t\sqrt\gamma)(t+\sqrt\gamma)^3},
\end{eqnarray}
and we obtain, after some calculation, 
\begin{equation}
  \Delta''(a_0) = - \frac{t^2(t+t\gamma+2\sqrt\gamma)^3\sqrt\gamma}
{2(1+t\sqrt\gamma)^2(t+\sqrt\gamma)^2}.
\end{equation}
By Taylor's formula, for $a\ge a_0$, 
\begin{equation}
  \Delta(a)= \frac12 \Delta''(a_0)(a-a_0)^2 + 
\frac16 \Delta'''(\widetilde{a}) (a-a_0)^3,
\end{equation}
for some $a_0\le \widetilde{a}\le a$.
For $\delta>0$ and $1\le \gamma_2<\infty$, set 
\begin{equation}
  C:= \sup \{ |\Delta'''(\widetilde{a},\gamma)| :
1\le \gamma\le \gamma_2, 
a_0(\gamma) \le \widetilde{a}\le a_0(\gamma)(1+\delta) \},
\end{equation}
where we have made the dependence on $\gamma$ of $\Delta'''$ 
explicit.
It follows from the smooth dependence of $\Delta(a,\gamma)$ 
on $\gamma$, as well as on $a$, that given $\gamma_2$, 
we can choose $\delta=\delta_2$ such that $C<\infty$.
Therefore
\begin{equation}
    \Delta(a)= \frac12 \Delta''(a_0)(a-a_0)^2 
\bigl(1+O(|a-a_0|\bigr),
\end{equation}
where $O(|a-a_0|)$ is uniform for $1\le \gamma\le \gamma_2$ 
and $a_0\le a\le a_0(1+\delta_2)$.
Indeed $O(|a-a_0|) \le C|a-a_0|$.
\end{proof}


In order to prove Proposition \ref{prop1}, we use \eqref{e-Yest}.
As 
\begin{equation}
  \cos\theta_c= \frac{x(a)^2-r(a)^2-t^{-2}}{2r(a)t^{-1}},
\end{equation}
and $\cos\theta_c\to -1$ as $a\to a_0$, 
we see that $\sin\frac{\theta_c}2=\sqrt{\frac{1-\cos\theta_c}2}$ 
is a smooth function of $a$ in  $[a_0,\infty)$.
Thus $\sin\frac{\theta_c}2= 1+O(|a-a_0|)$ for $a$ near $a_0$, 
$a\ge a_0$.
Inserting this information into \eqref{e-Yest} and using \eqref{e-Delest}, 
we obtain \eqref{e-Yexactest}.

\section{Proof of Theorem \ref{mainthm}}\label{sec-proof}

Take 
\begin{equation}
   L > \frac{2L_0}{a_0^{4/3}b_0},  
\qquad 0< \delta < \frac{\delta_0}{a_0b_0(1+\delta_0)},
\end{equation}
where $L_0, \delta_0$ are given in Proposition \ref{prop1}.
Let 
\begin{equation}
   n= \frac1{a_0} N  - x b_0 N^{1/3}.
\end{equation}
Set
\begin{equation}
  b:= \frac1{a_0}N - \frac{L_0}{a_0^{4/3}}N^{1/3}.
\end{equation}
Then for $L\le x\le \delta N^{2/3}$, 
the condition \eqref{e-kcond} for $a=\frac{N}{k}$ in Proposition \ref{prop1} 
is satisfied
for any $k$ in the range $n\le k\le b$, 
Following \cite{Merkl2}, we consider 
\begin{equation}
  \log \frac{\Prob(G([\gamma N],N)\le n)}{\Prob(G([\gamma N],N)\le b)}
\end{equation}
which equals 
\begin{equation}\label{e-sum}
  \sum_{k=n+1}^b \log (-Y_{21}(0;k))
\end{equation}
by \eqref{e-Gprod}.
Inserting \eqref{e-Yexactest} into \eqref{e-sum}, we obtain 
(cf. \cite{Merkl2}) 
\begin{equation}
\begin{split}
  \log \frac{\Prob(G([\gamma N],N)\le n)}{\Prob(G([\gamma N],N)\le b)} 
&=   -\frac13c_2a_0^3b_0^3  x^3 + O(x^4N^{-2/3}) + O(\log x) \\
&= -\frac1{12} x^3 + O(x^4N^{-2/3}) + O(log x) 
\end{split}
\end{equation}
as $c_2a_0^3b_0^3=\frac14$.
But by the result of Johansson \cite{kurtj:shape}, 
as $N\to\infty$, 
\begin{equation}
  \Prob(G([\gamma N],N)\le b) = F(-L_0a_0^{-4/3}b_0^{-1})
+ o(1),
\end{equation}
where $F(x)$ is the Tracy-Widom distribution. 
This proves Theorem \ref{mainthm}.

\bibliographystyle{plain}
\bibliography{paper10}

\end{document}